\documentclass[11pt]{article}%
\usepackage{amsthm,amsmath,amsfonts}
\usepackage{amssymb, url}
\usepackage{graphicx}
\usepackage{amsmath}
\usepackage{amsfonts}
\usepackage{amssymb}%
\theoremstyle{plain}
\newtheorem{theorem}{Theorem}
\newtheorem{lemma}{Lemma}
\newtheorem{proposition}{Proposition}
\newtheorem{corollary}{Corollary}

\theoremstyle{remark}
\newtheorem{remark}{Remark}

\newcommand{\E}{\mathbf{E}}

\newcommand{\dd}{\mathrm{d}}
\newcommand{\p}{\mathbf{P}}

\begin{document}

\title{Darling--Erd\H{o}s theorem for L\'evy processes at zero}
\author{P\'eter Kevei\thanks{ University of Szeged, Bolyai Institute, Aradi
v\'ertan\'uk tere 1, 6720 Szeged, Hungary, e-mail:
\texttt{kevei@math.u-szeged.hu}}
\and David M. Mason\thanks{ Department of Applied Economics and Statistics,
University of Delaware, 213 Townsend Hall, Newark, DE 19716, USA, e-mail:
\texttt{davidm@udel.edu}} }
\maketitle

\begin{abstract}
We establish two equivalent versions of the Darling--Erd\H{o}s theorem for
L\'evy processes in the domain of attraction of a stable process at zero
with index $\alpha\in(0,2)$. In the course of our proof we obtain a number of maximal and exponential inequalities for general L\'evy processes, which should be of separate interest.
\end{abstract}

\section{Introduction} \label{s1}

Let $\{ \xi_{k} \} _{k\geq1}$ be a sequence of independent mean zero and
variance one random variables and for each $n\geq1$ set $S_{n}=\xi_{1}%
+\dots+\xi_{n}$. Darling and {E}rd{\H{o}}s \cite{DA} proved that if the third
absolute moments of the $\{ \xi_{k} \} _{k\geq1}$ are uniformly bounded then
for all $x$, as $n\rightarrow\infty$,
\begin{equation}
\mathbf{P}\left(  A( n)  \max_{1\leq k\leq n}S_{k}/\sqrt
{k}- B( n)  \leq x\right)  \rightarrow\exp\left(  -\exp\left(
-x\right)  \right)  , \label{DA}%
\end{equation}
where we use the notation for $T>0$, $A( T) =\left(  2LLT\right)  ^{1/2}$ and
$B( T) =\left(  2LLT\right)  ^{1/2}+2^{-1}LLLT-2^{-1}L\left(  4\pi\right)  $,
with $LT=\log\left(  T\vee e\right)  $. Such a limiting distribution result is
now often called a Darling--Erd\H{o}s theorem. Einmahl \cite{E89} showed in
the i.i.d.~mean zero and variance one case that for (\ref{DA}) to hold it is
necessary and sufficient that
\[
LLt \, \mathbf{E}\left\{  \xi_{1}^{2} 
1\{  \left\vert \xi_{1}\right\vert \geq t\}  \right\}  \rightarrow0 \text{, as }t\rightarrow\infty\text{. }%
\]
Einmahl and Mason \cite{EM89} have obtained martingale Darling--Erd\H{o}s
theorems, and recently Dierickx and Einmahl \cite{DGE} have established
multivariate versions. Corresponding results for Brownian motion were established by Khoshnevisan et al.~\cite{Khos}. 

In the infinite-variance case Bertoin \cite{Bertoin98} proved 
Darling--Erd\H{o}s theorems for sums of i.i.d.~random variables from the normal domain of attraction of an $\alpha$-stable law. More precisely, if 
$\p ( \xi > x) \sim c x^{-\alpha}$ and $\p (\xi \leq -x ) = O(x^{-\alpha})$, as $x \to \infty$, for some $c > 0$ and $\alpha \in (0,1) \cup (1,2)$, and $\E \xi = 0$ for $\alpha > 1$ then for any $x \geq 0$
\begin{equation} \label{eq:Bertoin}
\lim_{n \to \infty} \p \left( \max_{k \leq n} k^{-1/\alpha} S_k \leq 
x (\log n )^{1/\alpha} \right) = e^{-c x^{-\alpha}}.
\end{equation}
Our work was motivated by this result.
In fact, our Theorem \ref{thm:DE} is a L\'{e}vy process version 
of Theorem 1 in \cite{Bertoin98}. \smallskip

Let $X_t, t\geq 0,$ be a L\'evy process in the domain of attraction of a stable process at zero with index $\alpha\in(0,2)$. Introduce the running supremum and the maximum jump process as
\begin{equation*}
\overline{X}_{t}=\sup_{s\leq t}X_{s},\quad m_{t}=\sup_{s\leq t}\Delta
X_{s}=\sup_{s\leq t}(X_{s}-X_{s-}).\label{eq:m-def}%
\end{equation*}
We consider for an appropriate positive increasing function $a(t)$ of $t>0$ the maximum of the scaled running supremum, the maximum of the scaled process, and the maximum of the scaled maximum jump process, defined as
\begin{equation}
Y_{t}=\sup_{t\leq s\leq1}\frac{\overline{X}_{s}}{a(s)},\quad 
Z_{t} =\sup_{t\leq s\leq1}\frac{X_{s}}{a(s)}, \quad 
M_{t}=\sup_{t\leq s\leq1} \frac{m_{s}}{a(s)}.\label{eq:M-Y-def}%
\end{equation}
For $\alpha=1$ the definitions of $Y$ and $Z$ are slightly different, see
Theorems \ref{thm:Y-lima1} and \ref{thm:DE}. Our goal is to derive analogues
of (\ref{DA}) and (\ref{eq:Bertoin}) for the L\'{e}vy process $X_{t}, t>0$. In particular,
we shall prove in our Theorem \ref{thm:Y-lim} that under suitable regularity conditions for all
$x>0$, in the case $\alpha\neq1$,
\[
\lim_{t\downarrow0}\mathbf{P}\left(  Y_{t}\,\left(  -\log{t}\right)
^{-1/\alpha}\leq x\right)  =e^{-x^{-\alpha}},
\]
and from this result we shall derive its Darling--Erd\H{o}s version in Theorem
\ref{thm:DE}
\[
\lim_{t\downarrow0}\mathbf{P}\left(  Z_{t}\,\left(  -\log{t}\right)
^{-1/\alpha}\leq x\right)  =e^{-x^{-\alpha}}.
\]
Along the way, in our Theorem \ref{thm:M-lim} we establish a similar result
for the scaled maximum jump process $M_{t}$. We fix our notation in Section
\ref{s2}, state our results in Section \ref{s3} and detail our proofs in
Sections \ref{s4} and \ref{s5}, where we derive some maximal and exponential inequalities for general L\'evy processes, which should of separate interest.

\section{Notation} \label{s2}

In this section we give our basic setup. Let $X_{t}$, $t\geq0$, be a L\'{e}vy
process with L\'{e}vy measure $\Lambda$ and without a normal component. Put
$\overline{\Lambda}_{+}(x)=\Lambda((x,\infty))$, $\overline{\Lambda}%
_{-}(x)=\Lambda((-\infty,-x))$, and for $u>0$ let
\begin{equation}
\varphi(u)=\sup\{x:\overline{\Lambda}_{+}(x)>u\}. \label{phi}%
\end{equation}
Note that $\overline{\Lambda}_{+}(x)>u$ iff $\varphi(u)>x$. Let $N$ be a
Poisson random measure on $(0,1)\times\mathbb{R}$ with intensity measure
$\mu(\mathrm{d}t,\mathrm{d}x)=\mathrm{d}t\times\Lambda(\mathrm{d}x)$ and let
$\widetilde{N}(\mathrm{d}t,\mathrm{d}y)=N(\mathrm{d}t,\mathrm{d}%
y)-\mathrm{d}t\Lambda(\mathrm{d}y)$ be the compensated Poisson measure. By the
L\'{e}vy-It\^{o} representation for suitable shift parameters $\gamma_{+}$ and
$\gamma_{-}$, with $\gamma=\gamma_{+}+\gamma_{-}$,
\begin{equation}%
\begin{split}
X_{t}  &  =\gamma t+\int_{0}^{t}\int_{|y|>1}yN(\mathrm{d}s,\mathrm{d}%
y)+\int_{0}^{t}\int_{|y|\leq1}y\widetilde{N}(\mathrm{d}s,\mathrm{d}y)\\
&  =\gamma_{+}t+\int_{0}^{t}\int_{(1,\infty)}yN(\mathrm{d}s,\mathrm{d}%
y)+\int_{0}^{t}\int_{(0,1]}y\widetilde{N}(\mathrm{d}s,\mathrm{d}y)\\
&  \quad+\gamma_{-}t+\int_{0}^{t}\int_{(-\infty,-1)}yN(\mathrm{d}%
s,\mathrm{d}y)+\int_{0}^{t}\int_{[-1,0)}y\widetilde{N}(\mathrm{d}%
s,\mathrm{d}y)\\
&  =:X_{t}^{+}+X_{t}^{-}.
\end{split}
\label{eq:X-repr}%
\end{equation}
We assume that $X^{+}$ belongs to the domain of attraction at zero of an
$\alpha$-stable law for some $\alpha\in(0,2)$, which means that~for some
norming and centering functions $a(t),c(t)$
\begin{equation}
\frac{X_{t}^{+}-c(t)}{a(t)}\overset{\mathcal{D}}{\longrightarrow}%
X,\quad\text{as }t\downarrow0, \label{eq:domain}%
\end{equation}
where $X$ is an $\alpha$-stable law. This happens if and only if
\begin{equation}
\overline{\Lambda}_{+}(x)=x^{-\alpha}\ell(x), \label{LL}%
\end{equation}
where $\ell$ is a slowly varying function at $0$; see Bertoin 
\cite[p.82]{Bertoin}, Maller and Mason \cite[Theorem 2.3]{MM08}. 
%Note that the spectrally positive and negative parts,
% $X_{t}^{+}$ and $X_{t}^{-}$, are independent. 
In what follows we assume that the constants $\gamma_{\pm}$ are
chosen such that
\begin{equation}
\gamma_{+}=%
\begin{cases}
\int_{(0,1]}y\Lambda(\mathrm{d}y), & \text{if } \int_{(0,1]}y\Lambda
(\mathrm{d}y) < \infty,\\
0, & \text{otherwise},
\end{cases}
\ \ \gamma_{-}=%
\begin{cases}
\int_{[-1,0)}y\Lambda(\mathrm{d}y), & \text{if } \int_{[-1,0)} |y|
\Lambda(\mathrm{d}y) < \infty,\\
0, & \text{otherwise.}%
\end{cases}
\label{eq:gamma}%
\end{equation}
Note that the integral $\int_{(0,1]} y \Lambda(\mathrm{d} y)$ is always finite
for $\alpha\in(0,1)$ and infinite for $\alpha\in(1,2)$, while for $\alpha= 1$
both cases can happen.

Without loss of generality we assume that $a(t)$ in (\ref{eq:domain}) is
increasing, moreover
\begin{equation}
a(t)=\varphi(1/t).
\label{eq:a-phi}%
\end{equation}
Using Remark (i) on page 320 of \cite{MM08} the function $c(t)$ in
(\ref{eq:domain}) can be chosen as
\begin{equation}
c(t)=t\nu(a(t))=t\left(  \gamma_{+}- \int_{(a(t),1]} y \Lambda(\mathrm{d}
y)\right)  , \label{ct}%
\end{equation}
where for $y>0$
\[
\nu(y)=\gamma_{+}-\int_{(y,1]} u \Lambda(\mathrm{d}u).
\]
For $0<\alpha<2$, with $\alpha\neq1$, it can be shown using standard
properties of regularly varying functions that, by the choice of $\gamma_{+}$,
\begin{equation*}\label{lima}
\lim_{t \downarrow0} \frac{c ( t)}{a( t)} =\frac{\alpha}{1-\alpha}.
\end{equation*}
This says that (\ref{eq:domain}) holds with $c( t) =0$ when $0<\alpha<2$, with
$\alpha\neq1$.

\section{Results} \label{s3}

From the monotonicity of $a$ it is simple that
\begin{equation*}
M_{t}=\sup_{t\leq s\leq1}\frac{m_{s}}{a(s)}=\sup\left\{  \frac{\Delta X_{s}%
}{a(s)}\,:s\in(t,1],\Delta X_{s}>0\right\}  \vee\frac{m_{t}}{a(t)},
\label{eq:M-form}%
\end{equation*}
where $a\vee b=\max\{a,b\}$. This simple observation allows us to calculate
the distribution of $M_{t}$. Indeed, for $x>0$ put
\[%
\begin{split}
A_{t,x}  &  =\left\{  (u,y):\frac{y}{a(u)}>x,\,u\in(t,1]\right\} \\
B_{t,x}  &  =\left\{  (u,y):\frac{y}{a(t)}>x,\,u\in(0,t]\right\}  .
\end{split}
\]
Then, recalling the definition of $N$ in (\ref{eq:X-repr}),
\[
\mathbf{P} \left(  \sup\left\{  \frac{\Delta X_{s}}{a(s)}\,: 
s \in (t,1],\Delta X_{s}>0\right\}  \leq x\right)  =\mathbf{P}(N(A_{t,x}%
)=0)=e^{-\mu(A_{t,x})},
\]
and
\[
\mathbf{P}(m_{t}\leq a(t)x)=\mathbf{P}(N(B_{t,x})=0)=e^{-\mu(B_{t,x})}.
\]
As $\mu(\mathrm{d}t,\mathrm{d}x)=\mathrm{d}t\times\Lambda(\mathrm{d}x)$, we
have
\[
\mu(A_{t,x})=\int_{t}^{1}\overline{\Lambda}_{+}(a(u)x)\mathrm{d}u\text{ and
}\mu(B_{t},x)=t\overline{\Lambda}_{+}(a(t)x).
\]
Since $A_{t,x}$ and $B_{t,x}$ are disjoint, we obtain
\begin{equation}
\mathbf{P}(M_{t}\leq x)=\exp\left\{  -\int_{t}^{1}\overline{\Lambda}%
_{+}(a(u)x)\mathrm{d}u-t\overline{\Lambda}_{+}(a(t)x)\right\}  .
\label{eq:M-df1}%
\end{equation}

\begin{remark}
If $X$ is a spectrally positive $\alpha$-stable process, $\alpha\in(0,2)$,
with $\overline{\Lambda}_{+}(x)=x^{-\alpha}$, then $\varphi(u)=u^{-1/\alpha}$,
$a(t)=t^{1/\alpha}$. Substituting into (\ref{eq:M-df1}) short calculation
gives
\[
\mathbf{P}(M_{t}\leq x)=\exp\left\{  -x^{-\alpha}(1 - \log t)\right\}  .
\]
Therefore, we obtain for any fixed $t>0$ the scaled maximum has Fr\'{e}chet
distribution, i.e.
\begin{equation}
\label{eq:stable-Fr}\mathbf{P}\left(  M_{t}\leq x ( 1- \log{t})^{1/\alpha}
\right)  = e^{-x^{-\alpha}}.
\end{equation}
\end{remark}

In what follows, we show that (\ref{eq:stable-Fr}) remains true in the limit
as $t \downarrow0$ for L\'evy processes in the domain of attraction of a
stable law at zero under regularity.

A measurable function $\ell$ is super-slowly varying at 0 with auxiliary function $\xi$,
if for some $\Delta>0$
\begin{equation}
\lim_{t\downarrow0}\frac{\ell(t\xi(t)^{\delta})}{\ell(t)}=1\quad
\text{uniformly in }\delta\in\lbrack0,\Delta]. \label{eq:ssv-def}%
\end{equation}
This is exactly the definition in Bingham et al.~\cite[Section 3.12.2]{BGT},
changing $x$ to $t^{-1}$ and $\xi(x)$ to $\xi(t^{-1})^{-1}$. See also
\cite[Section 2.3]{BGT}. We further assume that $\lim_{t\downarrow0}\xi(t)=0$
and that $\xi$ is nondecreasing in $\left(  0,c\right)  $ for some $c>0.$ If
(\ref{eq:ssv-def}) holds for some $\Delta>0$, and $\xi$ is nondecreasing then
(\ref{eq:ssv-def}) holds for any $\Delta>0$; see \cite[p.186]{BGT}.
In what follows we fix the function $\xi(t) = (-\log t)^{-1}$.

\begin{theorem}
\label{thm:M-lim} Assume that for $x>0$, $\overline{\Lambda}_{+}%
(x)=x^{-\alpha}\ell(x)$, $\alpha\in(0,2)$, where $\ell$ is a super-slowly
varying function at $0$ with auxiliary function $\xi(t)=(- \log t )^{-1}$.
Then for all $x>0$
\begin{equation}
\lim_{t\downarrow0}\mathbf{P}\left(  M_{t}\,\left(  -\log{t}\right)
^{-1/\alpha}\leq x\right)  =e^{-x^{-\alpha}}. \label{ED}%
\end{equation}
\end{theorem}

\begin{remark}
The super-slowly varying condition is not very restrictive. The slowly varying
functions $\ell(t) = (-\log t)^{\beta}$, $\beta> 0$, $\ell(t) = \exp\{ (-\log
t)^{\beta} \}$, $\beta\in(0,1)$ are super-slowly varying with auxiliary
function $\xi(t) = (- \log t)^{-1}$. The function $\ell(t) = \exp\{ (-\log t)
/ \log(- \log t) \}$ is slowly varying, but not super-slowly varying with
auxiliary function $\xi$.
\end{remark}

\begin{remark}
We also note that Theorem \ref{thm:M-lim} is a result on the maximum of
a Poisson point process, therefore $\Lambda(\mathrm{d}x)$ does not have to be a
L\'{e}vy measure. Thus Theorem \ref{thm:M-lim} remains true for any $\alpha>0$.
\end{remark}

For our next result assume that the spectrally negative part does not dominate
in the sense
\begin{equation}
\limsup_{x\downarrow0}\frac{\overline{\Lambda}_{-}(x)}{\overline{\Lambda}%
_{+}(x)}<\infty. \label{eq:lambda+-}%
\end{equation}

\begin{theorem}
\label{thm:Y-lim} Assume that $X_{t}$ is a L\'{e}vy process without normal
component such that for $x>0$, $\overline{\Lambda}_{+}(x)=x^{-\alpha}\ell(x)$,
$\alpha\in(0,2)$ with $\alpha\neq1$, where $\ell$ is a super-slowly varying
function at 0 with auxiliary function $\xi(t)=(- \log t )^{-1}$, and
(\ref{eq:lambda+-}) holds. Then for all $x>0$
\[
\lim_{t\downarrow0}\mathbf{P}\left(  Y_{t}\,\left(  -\log{t}\right)
^{-1/\alpha}\leq x\right)  =e^{-x^{-\alpha}}.
\]
\end{theorem}

This result also holds for $\alpha= 1$ but, as usual, a different centering is
needed. As in (\ref{ct}), for $\alpha= 1$ let
\begin{equation}
\label{eq:c-a1}c(t) =
\begin{cases}
t \int_{(0,{a(t)}]} y \Lambda ( \mathrm{d} y ), & \text{if } \ \int_{(0,1]}
y \Lambda(\mathrm{d} y) < \infty,\\
- t \int_{(a(t),1]} y \Lambda ( \mathrm{d} y), & \text{if } \ \int_{(0,1]}
y \Lambda(\mathrm{d} y) = \infty.
\end{cases}
\end{equation}

\begin{theorem}
\label{thm:Y-lima1} Assume that $X_{t}$ is a L\'{e}vy process without normal
component such that for $x>0$, $\overline{\Lambda}_{+}(x)=x^{-1}\ell(x)$,
where $\ell$ is a super-slowly varying function at 0 with auxiliary function
$\xi(t)=(- \log t )^{-1}$, (\ref{eq:lambda+-}) holds, and $\int_{[-1,0)} -y
\Lambda(\mathrm{d} y) < \infty$. Then for all $x>0$
\[
\lim_{t\downarrow0} \mathbf{P}\left(  \sup_{t \leq s \leq1} \frac{ \sup_{u
\leq s} ( X_{u} - c(u)) }{a(s)} \,\left(  -\log{t}\right)^{-1} \leq
x\right)  =e^{-x^{-1}}.
\]
\end{theorem}

As a consequence, we obtain the following Darling--Erd\H{o}s result.

\begin{theorem}
\label{thm:DE} Assume that $X_{t}$ is a L\'{e}vy process without normal
component such that for $x>0$, $\overline{\Lambda}_{+}(x)=x^{-\alpha}\ell(x)$,
$\alpha\in(0,2)$, where $\ell$ is a super-slowly varying function at 0 with
auxiliary function $\xi(t)=(-\log t)^{-1}$, and (\ref{eq:lambda+-}) holds. For
$\alpha=1$ additionally assume $\int_{[-1,0)}-y\Lambda(\mathrm{d}y)<\infty$.
Then for all $x>0$
\begin{equation} \label{eq:DElim}
\lim_{t\downarrow0}\mathbf{P}\left(  \sup_{t\leq s\leq1}\frac{X_{s}%
-c(s)}{a(s)}\,\left(  -\log{t}\right)  ^{-1/\alpha}\leq x\right)
=e^{-x^{-\alpha}},
\end{equation}
where $c(s)\equiv0$ for $\alpha\neq1$, and given in (\ref{eq:c-a1}) for
$\alpha=1$.
\end{theorem}

\begin{remark}
We note that the conditions for the corresponding result for sums of i.i.d.~random variables in \cite[Theorem 1]{Bertoin98} are more stringent. The non-dominating negative tail assumption is the same as (\ref{eq:lambda+-}), but in \cite{Bertoin98} it is assumed that the slowly varying function $\ell$ in (\ref{LL}) is constant, and the $\alpha = 1$ case is excluded. It will be apparent from the proofs that the nontrivial slowly varying function significantly complicates the arguments.

We also mention that large time results similar to (\ref{eq:DElim}) for stable processes are stated in Theorem 5 of \cite{Bertoin98} based on the correspondence between stable processes and stable Ornstein--Uhlenbeck processes. 
Theorem 5 in \cite{Bertoin98} can be deduced from Corollary 5.3 in Rootz\'en \cite{Rootzen}, since Ornstein--Uhlenbeck processes can be represented as stable moving average processes with exponential kernel function (see e.g.~Applebaum \cite[Section 4.3.5]{Applebaum}).
\end{remark}

\section{Proof of Theorem \ref{thm:M-lim}} \label{s4}

Since here the spectrally negative part does not play a role, to ease the
notation we suppress the lower index, i.e.~$\overline{\Lambda}=\overline
{\Lambda}_{+}$. From (\ref{eq:M-df1}) we get for fixed $x>0$
\begin{equation}%
\begin{split}
\mathbf{P}\left(  M_{t}\left(  -\log{t}\right)  ^{-1/\alpha}\leq x\right)
=\exp\left\{  -\int_{t}^{1}\overline{\Lambda}(a(u)(- \log t)^{1/\alpha
}x)\mathrm{d}u- t \overline{\Lambda}(a(t)(-\log t)^{1/\alpha}x)\right\}  .
\end{split}
\label{H}%
\end{equation}
In what follows, we need that
\begin{equation}
\overline{\Lambda} (a(u))\sim u^{-1}\text{, as } u\downarrow0.
\label{eq:a-prop1}%
\end{equation}
To see this, define for $x>0$
\[
f\left(  x\right)  =\overline{\Lambda} (1/x)=x^{\alpha}\ell(1/x).
\]
Clearly $f$ is increasing and regularly varying with index $\alpha$ at
$\infty$. Recall (\ref{phi}) and  set for $y>0$
\[%
\begin{split}
f^{-1}\left(  y\right)   &  =\inf\left\{  x:f\left(  x\right)  >y\right\} \\
&  =\inf\left\{  x:\overline{\Lambda} (1/x)>y\right\} \\
&  =1/\sup\left\{  x^{-1}:\overline{\Lambda} (1/x)>y\right\}
=1/\varphi(y).
\end{split}
\]
By (\ref{eq:a-phi})
and  Theorem 1.5.12 of \cite{BGT} we have that as $y\rightarrow\infty$
\[
f\left(  f^{-1}\left(  y\right)  \right)  =
\overline{\Lambda} (\varphi(y))\sim y,
\]
which by the change of variable $y=u^{-1}$ gives (\ref{eq:a-prop1}).

Let $h$ be an auxiliary function to be chosen later, which is continuous,
increasing on $\left(  0,1\right)  $ and $1>h(t)>t$. We can write the exponent
in (\ref{H}) as%
\begin{equation}%
\begin{split}
&  -\int_{t}^{h(t)}\overline{\Lambda}(a(u)(- \log t)^{1/\alpha}x)\mathrm{d}%
u-\int_{h(t)}^{1}\overline{\Lambda}(a(u)( - \log t)^{1/\alpha}x)\mathrm{d}u\\
&  -t\overline{\Lambda}(a(t)(- \log t )^{1/\alpha}x).
\end{split}
\label{EH}%
\end{equation}
By the assumption on $\overline{\Lambda}$
\begin{equation}
\frac{\overline{\Lambda}\left(  a(u)x(- \log t )^{1/\alpha}\right)
}{\overline{\Lambda}\left(  a(u)\right)  }=x^{-\alpha}(- \log t )^{-1}%
\frac{\ell\left(  a(u)x(- \log t )^{1/\alpha}\right)  }{\ell(a(u))}.
\label{eq:intLambda-1}%
\end{equation}
By the definition of super-slowly varying functions, for any $\varepsilon>0$
there exists $t_{0}>0$ such that for all $s\in(0,t_{0})$
\begin{equation}
\sup_{\xi(s)^{\Delta}\leq y\leq1}\left\vert \frac{\ell(sy)}{\ell
(s)}-1\right\vert <\varepsilon, \label{eq:ssv-aux}%
\end{equation}
where $\Delta> 0$. To see this note that for any $\xi(s)^{\Delta}\leq y\leq1$
there exists a $0\leq\rho\leq\Delta$ such that $\xi(s)^{\rho}=y$. We choose
for $0<\beta<1$ and $t\in(0,1)$%
\begin{equation*}
h\left(  t\right)  =h_{\beta}(t)=\exp\left\{  -\left(  - \log t \right)
^{\beta}\right\}  . \label{eq:h}%
\end{equation*}
We claim that
\begin{equation}
\lim_{t\downarrow0}\sup_{u\in\lbrack t,h(t)]}\left\vert \frac{\ell(a(u))}%
{\ell\left(  a(u)x(- \log t )^{1/\alpha}\right)  }-1\right\vert =0.
\label{eq:ssv-claim}%
\end{equation}
In (\ref{eq:ssv-aux}) choose
\begin{equation}
s=s(u,t,x)=a(u)x\left(  \log t^{-1}\right)  ^{1/\alpha}\quad\text{and
}\ y=y(t,x)=x^{-1}\left(  \log t^{-1}\right)  ^{-1/\alpha}. \label{eq:sy-def}%
\end{equation}
Clearly, $y\leq1$ for $t$ small enough. Thus, in order to use
(\ref{eq:ssv-aux}) we have to check that
\begin{equation}
\lim_{t\downarrow0}\sup_{u\in\lbrack t,h(t)]}a(u)x\left(  \log t^{-1}\right)
^{1/\alpha}=0, \label{eq:ssv-cond1}%
\end{equation}
and, with $s,y$ in (\ref{eq:sy-def}) and $u\in\lbrack t,h(t)]$,
\begin{equation}
\xi(s)^{\Delta}\leq y=x^{-1}\left(  \log t^{-1}\right)  ^{-1/\alpha}.
\label{eq:ssv-cond2}%
\end{equation}
Since $a$ is regularly varying at $0$ with parameter $1/\alpha$
\begin{equation}
\log\left(  a(h(t))(- \log t )^{1/\alpha}\right)  \sim\frac{-\left(  - \log t
\right)  ^{\beta}}{\alpha},\quad\text{as }t\downarrow0. \label{eq:h-aux2}%
\end{equation}
Using the monotonicity of $a$ and (\ref{eq:h-aux2}), for $u\in[t,h(t)]$, $t$
small enough
\begin{equation}
a(u)(- \log t )^{1/\alpha}\leq a(h(t))(- \log t )^{1/\alpha}\leq\exp\left\{
\frac{-\left(  - \log t \right)  ^{\beta}}{2\alpha}\right\}  .
\label{eq:h-cond1}%
\end{equation}
The latter upper bound tends to 0 as $t\downarrow0$, therefore
(\ref{eq:ssv-cond1}) follows.

By (\ref{eq:h-cond1}) for any $\Delta>1$ and $t>0$ small enough
\begin{equation}
\xi\left(  a(h(t))x(- \log t )^{1/\alpha}\right)  ^{\Delta}\leq{\left(
\frac{1}{3\alpha}\log\left(  1/h(t)\right)  \right)  ^{-\Delta}}%
=(3\alpha)^{\Delta}(- \log t )^{-\beta\Delta}. \label{eq:xi-ineq}%
\end{equation}
By the monotonicity of $\xi$ and $a$, and by (\ref{eq:xi-ineq}) we have
\[%
\begin{split}
\sup_{t\leq u\leq h(t)}\xi(s)^{\Delta}  &  =\sup_{t\leq u\leq h(t)}\xi\left(
a(u)x(\log(1/t))^{1/\alpha}\right)  ^{\Delta}=\xi\left(  a(h(t))x(\log
(1/t))^{1/\alpha}\right)  ^{\Delta}\\
&  \leq(3\alpha)^{\Delta}(- \log t )^{-\beta\Delta}\leq x^{-1}(\log
1/t)^{-1/\alpha},
\end{split}
\]
where the last inequality holds for $t$ large enough if $\beta\Delta>1/\alpha
$. Since, by the remark before Theorem \ref{thm:M-lim}, $\Delta$ can be chosen
to be large, (\ref{eq:ssv-cond2}) holds, and (\ref{eq:ssv-claim}) follows.

Thus, by (\ref{eq:intLambda-1}) uniformly in $u\in(t,h(t))$
\[
\frac{\overline{\Lambda}\left(  a(u)x(- \log t )^{1/\alpha}\right)
}{\overline{\Lambda}\left(  a(u)\right)  }\sim x^{-\alpha}(- \log t )^{-1}.
\]
Therefore, using also (\ref{eq:a-prop1}),
\[%
\begin{split}
\int_{t}^{h(t)}\overline{\Lambda}\left(  a(u)x(- \log t )^{1/\alpha}\right)
\mathrm{d}u  &  \sim x^{-\alpha}(- \log t )^{-1}\int_{t}^{h(t)}u^{-1}%
\mathrm{d}u\\
&  \sim x^{-\alpha}\left(  1-\frac{\log1/h(t)}{- \log t }\right)  \sim
x^{-\alpha},\text{ as }t\downarrow0.
\end{split}
\]
Next we see that%
\[
\int_{h(t)}^{1}\overline{\Lambda}\left(  a(u)x(- \log t )^{1/\alpha}\right)
\mathrm{d}u=\int_{h(t)}^{1}\overline{\Lambda}\left(  a(u)\right)  \left[
\frac{\overline{\Lambda}\left(  a(u)x(- \log t )^{1/\alpha}\right)
}{\overline{\Lambda}\left(  a(u)\right)  }\right]  \mathrm{d}u,
\]
which by (\ref{eq:intLambda-1})
\[
=x^{-\alpha}(- \log t )^{-1}\int_{h(t)}^{1}\overline{\Lambda}\left(
a(u)\right)  \left[  \frac{\ell\left(  a(u)x(- \log t )^{1/\alpha}\right)
}{\ell\left(  a(u)\right)  }\right]  \mathrm{d}u.
\]
Applying part (ii) of Theorem 1.5.6 in \cite{BGT} we see that this last bound
is for any $\delta>0$, some $A_{\delta}>0$ and for all small enough $t>0$%
\[
\leq A_{\delta}x^{-\alpha}(- \log t )^{-1}\int_{h(t)}^{1}\overline{\Lambda
}\left(  a(u)\right)  \left(  x(- \log t )^{1/\alpha}\right)  ^{\delta
}\mathrm{d}u.
\]
By (\ref{eq:a-prop1}) we can infer that there exists a $B>0$ such that for all
$u\in(0,1]$%
\[
\overline{\Lambda}(a(u))\leq Bu^{-1}\text{.}%
\]
Thus with $C=A_{\delta}B$
\[%
\begin{split}
&  A_{\delta}x^{-\alpha}(- \log t )^{-1}\int_{h(t)}^{1}\overline{\Lambda
}\left(  a(u)\right)  \left(  x(- \log t )^{1/\alpha}\right)  ^{\delta
}\mathrm{d}u\\
&  \leq Cx^{\delta-\alpha}(- \log t )^{\delta/\alpha-1}\int_{h_{\beta}(t)}%
^{1}u^{-1}\mathrm{d}u\\
&  =Cx^{\delta-\alpha}(- \log t )^{\delta/\alpha-1+\beta},
\end{split}
\]
which for small enough $\delta>0$ converges to zero as $t\downarrow0$.

Therefore it follows that
\[
\lim_{t\downarrow0}\int_{t}^{1}\overline{\Lambda}\left(  a(u)x (- \log
t)^{1/\alpha}\right)  \mathrm{d}u=x^{-\alpha}.
\]
Finally, the for third term in (\ref{EH}) we have, by (\ref{eq:a-prop1}) and
(\ref{eq:ssv-claim}),%

\[
t\overline{\Lambda}\left(  a(t)x(- \log t )^{1/\alpha}\right)  \sim
\frac{\overline{\Lambda}\left(  a(t)x(- \log t )^{1/\alpha}\right)
}{\overline{\Lambda}\left(  a(t)\right)  }\sim x^{-\alpha}\left(  - \log t
\right)  ^{-1},
\]
which converges to zero as $t\downarrow0$, and statement (\ref{ED}) follows.
%\hspace*{1pt} \hfill\qed

\section{Proofs of Theorems \ref{thm:Y-lim}, \ref{thm:Y-lima1}, and
\ref{thm:DE}}

\label{s5}

\subsection{Exponential inequalities for general L\'evy processes}

\noindent In this subsection for convenience of presentation we state and
prove the exponential inequalities that are needed in the proof of Theorem
\ref{thm:Y-lim}. All of them are derived from Proposition \ref{prop} below,
which may be of separate interest.\medskip

Let $X_{t}$, $t\geq0$, be a L\'{e}vy process without a normal component with
L\'{e}vy measure $\Lambda$. As before for $x>0$, $\overline{\Lambda}%
_{+}(x)=\Lambda((x,\infty))$ and $\overline{\Lambda}_{-}(x)=\Lambda
((-\infty,-x))$. For any fixed $a>0$ introduce the L\'{e}vy processes
\begin{equation}
X_{s}^{(a)}=\int_{0}^{s}\int_{(0,a]}y\widetilde{N}(\mathrm{d}u,\mathrm{d}y)
\ \text{ and } \ X_{s}^{(-a)}=\int_{0}^{s}\int_{\left[  -a,0\right)
}y\widetilde{N}(\mathrm{d}u,\mathrm{d}y), \quad s\geq0. \label{NA}%
\end{equation}
%where we assume both $\overline{\Lambda}_{+}(0+)=\infty$ and 
%$\overline {\Lambda}_{-}(0)=\infty$.\smallskip

Set for $a\geq0$%
\begin{equation}
\label{eq:def_B}B( a) =\int_{0}^{a}y^{2}\Lambda(\mathrm{d}y) \ \text{ and }
\ B( -a) =\int_{-a}^{0}y^{2}\Lambda(\mathrm{d}y).
\end{equation}
We note that the following proposition holds for general L\'evy process,
regular variation of the L\'evy measure is not needed here.

\begin{proposition}
\label{prop} For all $a>0$, $b>0$, $p\geq1$ integer, and $0 < t$%
\begin{equation}
\mathbf{P} \left(  \sup_{s\leq t}X_{s}^{(a)} > b\right)  
\leq \exp \left\{
\frac{b}{\left(  1+1/p\right)  a}\left(  1+\log\left(  
\frac{t B( a)  \left(  p!\right)  ^{1/p}}{ab}\right)  \right)  \right\}
\label{eq1}%
\end{equation}
and for all $a>0$, $b>0$ and $0 < t$
\begin{equation}
\mathbf{P}\left(  \sup_{s\leq t}\left(  -X_{s}^{(a)}\right)  >b\right)
\leq\exp\left(  -\frac{b^{2}}{2 t B ( a) }\right)  . \label{eq2}%
\end{equation}
Moreover, inequality (\ref{eq1}) holds with $\sup_{s\leq t}X_{s}^{(a)}$
replaced by $\sup_{s\leq t}\left(  -X_{s}^{(-a)}\right)  $ and inequality
(\ref{eq2}) remains true with $\sup_{s\leq t}\left(  -X_{s}^{(a)}\right)  $
replaced by $\sup_{s\leq t}X_{s}^{(-a)}$, and where $B( a)$ is
replaced by $B(-a)$ in both cases.
\end{proposition}

\noindent\textit{Proof.} We shall borrow steps from the proof of Lemma 1 of
Sato \cite{SATO}. Clearly, $X_{s}^{(a)}$ is a martingale, thus by Doob's
martingale inequality, for any $\theta>0$
\begin{equation}%
\begin{split}
\mathbf{P}\left(  \sup_{s\leq t}X_{s}^{(a)}>b\right)   &  =\mathbf{P}\left(
\exp\left\{  \theta\sup_{s\leq t}X_{s}^{(a)}\right\}  >e^{\theta b}\right) \\
&  \leq e^{-\theta b}\mathbf{E}\exp\{\theta X_{t}^{(a)}\}.
\end{split}
\label{eq:doob1}%
\end{equation}
The difficult issue here is to choose the right $\theta$. \medskip

Set for $\theta\in\mathbb{R}$,%
\[
\xi_{t}\left(  \theta\right)  :=\log\mathbf{E}\exp\{\theta X_{t}^{(a)}%
\}=t\int_{0}^{a}\left(  e^{\theta y}-1-\theta y\right)  \Lambda(\mathrm{d}y).
\]
Since $\left\vert e^{v}-1-v\right\vert \leq v^{2}\exp\left(  \left\vert
v\right\vert \right)  /2$ for all $v\in\mathbb{R}$, we see that
\[
\left\vert \xi_{t} ( \theta) \right\vert \leq\frac{\theta^{2}}{2}\int_{0}%
^{a}y^{2}t\Lambda(\mathrm{d}y)\exp\left(  a\left\vert \theta\right\vert
\right)  <\infty, \quad \theta\in\mathbb{R}.
\]
Thus $\mathbf{E}\exp\{\theta X_{t}^{(a)}\}<\infty$ for all $\theta
\in\mathbb{R}$. Differentiating $\xi_{t} ( \theta) $ with respect to $\theta$
we obtain for all $\theta\in\mathbb{R}$
\[
\xi_{t}^{\prime} (\theta) =\int_{0}^{a}y\left(  e^{\theta y}-1\right)
t\Lambda(\mathrm{d}y),
\]
and differentiating again, for all $\theta\in\mathbb{R}$%
\begin{equation}
\xi_{t}^{{\prime\prime}} ( \theta) =\int_{0}^{a}y^{2}e^{\theta y}%
t\Lambda(\mathrm{d}y)>0, \label{incr}%
\end{equation}
from which we see that
\[
\xi_{t}^{\prime}\left(  \theta\right)  \downarrow-\mu=:-\int_{0}^{a}%
yt\Lambda(\mathrm{d} y)\text{, as } \ \theta\downarrow-\infty,
\]
where $-\infty\leq-\mu<0$, $\xi_{t}^{\prime}\left(  0\right)  =0$ and $\xi
_{t}^{\prime}\left(  \theta\right)  \uparrow\infty$, as $\theta\uparrow\infty
$.\medskip

For any $-\mu<x<\infty$ introduce the inverse to $\xi_{t}^{\prime}\left(
\theta\right)  :$%
\begin{equation}
\xi_{t}^{\prime}\left(  \eta_{t}\left(  x\right)  \right)  =x. \label{inver}%
\end{equation}
The function $\eta_{t}$ is well defined on $\left(  -\mu,\infty\right)  $,
since by (\ref{incr}), $\xi_{t}^{\prime}\left(  \theta\right)  $ is strictly
increasing and continuous as a function of $\theta$. Furthermore by the
inverse function theorem we have
\begin{equation}
\xi_{t}^{{\prime\prime}}( \eta_{t}( x) ) \eta_{t}^{\prime}( x) =1 \text{, for
}-\mu<x<\infty, \label{one}%
\end{equation}
and we know from the above that $\eta_{t}( x) >0$ if and only if $x>0$. Now by
(\ref{eq:doob1}) with $\theta=\eta_{t}( b) $ and (\ref{inver}) for any $b>0$,
\[
\begin{split}
\mathbf{P}\left(  \sup_{s\leq t}X_{s}^{(a)}>b\right)   &  \leq\mathbf{P}%
\left(  \eta_{t} ( b) \sup_{s\leq t}X_{s}^{(a)}>\eta_{t}( b) b \right) \\
&  \leq\exp\left\{  \xi_{t} ( \eta_{t}( b) ) - \eta_{t} ( b) \xi_{t}^{\prime}(
\eta_{t}( b) ) \right\}  .
\end{split}
\]
Observe that
\[
\begin{split}
&  \xi_{t} \left(  \eta_{t}\left(  b\right)  \right)  -\eta_{t}\left(
b\right)  \xi_{t}^{\prime}\left(  \eta_{t}\left(  b\right)  \right)  =\int
_{0}^{\eta_{t}\left(  b\right)  }\xi_{t}^{\prime}\left(  s\right)
\mathrm{d}s-\eta_{t}\left(  b\right)  \xi_{t}^{\prime}\left(  \eta_{t}\left(
b\right)  \right) \\
&  =-\int_{0}^{\eta_{t}\left(  b\right)  }s\xi_{t}^{^{\prime\prime}}\left(
s\right)  \mathrm{d}s=-\int_{0}^{b}\eta_{t}\left(  x\right)  \xi_{t}%
^{^{\prime\prime}}\left(  \eta_{t}\left(  x\right)  \right)  \eta_{t}^{\prime
}\left(  x\right)  \mathrm{d}x,
\end{split}
\]
which by (\ref{one}) is equal to%
\[
=-\int_{0}^{b}\eta_{t}\left(  x\right)  \mathrm{d}x.
\]
Thus for all $b>0$%
\begin{equation}
\mathbf{P}\left(  \sup_{s\leq t}X_{s}^{(a)}>b\right)  \leq\exp\left(
-\int_{0}^{b}\eta_{t}\left(  x\right)  \mathrm{d}x\right)  . \label{eta}%
\end{equation}
Since $\exp\left(  v\right)  -1\leq v\exp\left(  v\right)  $ for all $v\geq0$,
for $\tau\geq0$
\[
\begin{split}
\xi_{t}^{\prime}\left(  \tau\right)   &  =\int_{0}^{a}y\left(  e^{\tau
y}-1\right)  t\Lambda(\mathrm{d}y) \leq t\int_{0}^{a}y^{2}\Lambda(\mathrm{d}y)
\tau\exp\left(  a\tau\right) \\
&  = tB ( a )  \tau \exp \left(  a\tau \right)  \text{, }%
\end{split}
\]
from which it follows by (\ref{inver}) that
\[
x\leq t B ( a )  \eta_{t} ( x)  \exp\left(  a\eta_{t} ( x)  \right)
  \text{, for } \ x\geq0.
\]
The inequality $v\leq\left(  p!\right)  ^{1/p}\exp\left(  \frac{v}{p}\right)
$ for $p\geq1$ and $v\geq0$, gives
\[
x\leq\frac{tB\left(  a\right)  }{a}\left(  p!\right)  ^{1/p}\exp\left(
\left(  1+1/p\right)  a\eta_{t}\left(  x\right)  \right)  \text{, for
}x>0\text{,}%
\]
and thus
\[
\frac{\log x}{\left(  1+1/p\right)  a}-\frac{1}{\left(  1+1/p\right)  a}%
\log\left(  \frac{tB\left(  a\right)  \left(  p!\right)  ^{1/p}}{a}\right)
\leq\eta_{t}\left(  x\right)  \text{, for }x>0.
\]
Hence after a little algebra we get
\begin{equation*}%
\begin{split}
\exp\left(  -\int_{0}^{b}\eta_{t}\left(  x\right)  \mathrm{d}x\right)   &
\leq\exp\left(  -\int_{0}^{b} \frac{ \log x \ \mathrm{d}x}{\left(
1+1/p\right)  a} +\frac{b}{\left(  1+1/p\right)  a}\log\left(  \frac{tB\left(
a\right)  \left(  p!\right)  ^{1/p}}{a}\right)  \right) \\
&  =\exp\left\{  \frac{b}{\left(  1+1/p\right)  a}\left(  1+\log\left(
\frac{tB\left(  a\right)  \left(  p!\right)  ^{1/p}}{ab}\right)  \right)
\right\}  ,\label{inq}%
\end{split}
\end{equation*}
which on account of (\ref{eta}) gives (\ref{eq1}).\medskip

Next consider inequality (\ref{eq2}). The process $-X_{s}^{(a)}$, $s\geq0$, is
also a martingale. Therefore exactly as above for all $\theta>0$
\[
\mathbf{P}\left(  \sup_{s\leq t}\left(  -X_{s}^{(a)}\right)  >b\right)  \leq
e^{-\theta b}\mathbf{E}\exp\{-\theta X_{t}^{(a)}\}=\exp\left(  -\theta
b+\gamma_{t}\left(  \theta\right)  \right)  ,
\]
where $\gamma_{t}\left(  \theta\right)  =\xi_{t}\left(  -\theta\right)  $. We
get
\[
\gamma_{t}^{\prime}\left(  \theta\right)  =-\xi_{t}^{\prime}\left(
-\theta\right)  =t\int_{0}^{a}y\left(  1-e^{-\theta y}\right)  \Lambda
(\mathrm{d}y),
\]
and $\gamma_{t}^{\prime\prime}\left(  \theta\right)  =\xi_{t}^{^{\prime\prime
}}\left(  -\theta\right)  >0$, from which we see that
\[
\gamma_{t}^{\prime}\left(  \theta\right)  \uparrow\mu=t\int_{0}^{a}%
y\Lambda(\mathrm{d}y)\text{, as }\theta\uparrow\infty,
\]
where $0<\mu\leq\infty$, $\gamma_{t}^{\prime}\left(  0\right)  =0$ and
$\gamma_{t}^{\prime}\left(  \theta\right)  \downarrow-\infty$, as
$\theta\downarrow-\infty$. For any $-\infty<x<\mu$ introduce the inverse to
$\gamma_{t}^{\prime}\left(  \theta\right)  :$%
\begin{equation*}
\gamma_{t}^{\prime}\left(  \kappa_{t}\left(  x\right)  \right)  =x.
\label{inver1}%
\end{equation*}
The function $\kappa_{t}$ is well defined on $\left(  -\infty,\mu\right)  $,
since by $\gamma_{t}^{\prime\prime}\left(  \theta\right)  =\xi_{t}%
^{^{\prime\prime}}\left(  -\theta\right)  >0$, $\gamma_{t}^{\prime}\left(
\theta\right)  $ is strictly increasing and continuous as a function of
$\theta$. Furthermore by the inverse function theorem we have
\begin{equation*}
\gamma_{t}^{{\prime\prime}}\left(  \kappa_{t}\left(  x\right)  \right)
\kappa_{t}^{\prime}\left(  x\right)  =1\text{, for }-\infty<x<\mu,
\label{one1}%
\end{equation*}
and we know from the above that $\kappa_{t}\left(  x\right)  >0$ if and only
if $x>0$.\medskip

Now just as in the proof (\ref{eq1}), for all $b>0$
\begin{equation}
\mathbf{P}\left(  \sup_{0\leq s\leq t}\left(  -X_{s}^{(a)}\right)  >b\right)
\leq\exp\left(  -\int_{0}^{b}\kappa_{t}\left(  x\right)  \mathrm{d}x\right)  .
\label{kappa}%
\end{equation}
Since $1-\exp\left(  -v\right)  \leq v$ for $v>0$, we have for all $\theta
\geq0$,%
\begin{equation}%
\begin{split}
\gamma_{t}^{\prime}\left(  \theta\right)   &  =-\xi_{t}^{\prime}\left(
-\theta\right)  =t\int_{0}^{a}y\left(  1-e^{-\theta y}\right)  \Lambda
(\mathrm{d}y)\\
&  \leq t\theta\int_{0}^{a}y^{2}\Lambda(\mathrm{d}y)=\theta tB\left(
a\right)  \text{,}\label{kk}%
\end{split}
\end{equation}
from which it follows by setting $\theta=\kappa_{t}\left(  x\right)  $ into
(\ref{kk}) that $x\leq tB\left(  a\right)  \kappa_{t}\left(  x\right)  $. This
gives (\ref{eq2}) by (\ref{kappa}).\smallskip

The validity of the moreover part of the statement of Proposition \ref{prop}
is obvious. \qed 

\subsection{Applications of Proposition \ref{prop}}

In what follows we assume (\ref{LL}). Then Karamata's theorem implies that for
$B$ in (\ref{eq:def_B})
\begin{equation}
\label{eq:B-asy}B(y) = \int_{(0,y]} y^{2} \Lambda(\mathrm{d} y) \sim
\frac{\alpha}{2- \alpha} y^{2} \overline\Lambda_{+}(y) \quad\text{as } y
\downarrow0.
\end{equation}
For any $0<\beta<\alpha$, select $0<\rho<\rho_{2}$ small and $\kappa>1$
depending on $\alpha$ and $\beta$ so that for all $0<xu\leq\rho$ with $x\geq1$
such that by the Potter bounds [Theorem 1.5.6 \cite[Section 2.3]{BGT}, p. 25],%
\begin{equation}
\frac{\ell( ux) }{\ell( u ) } \leq\kappa x^{\beta}. \label{beta}%
\end{equation}

\begin{corollary}
\label{corr:1} Assume (\ref{LL}). For any $\varepsilon\in(0,1)$ there exist
$\alpha^{\prime}> \alpha$, $t_{0} > 0$, and $A > 0$, such that whenever
$\max\{ a(t) x, a(t), t \} < t_{0}$ and $x ( 1 - \varepsilon) > 1$
\begin{equation}
\mathbf{P}\left(  \sup_{s\leq t}X_{s}^{(a(t) x (1-\varepsilon))} > a(t) x ( 1
- \varepsilon/ 2) \right)  \leq A x^{-\alpha^{\prime}}. \label{prime}%
\end{equation}
\end{corollary}

{\noindent\textit{Proof.}} Let $a=a(t)x(1-\varepsilon)$ with $x(1-\varepsilon
)\geq1$, $b=a\left(  t\right)  x(1- {\varepsilon/2})$, and set $q=(1-
{\varepsilon}/{2})/\left(  1-\varepsilon\right) $. Choose the integer $p$ so
large that $1+1/p<q$. By (\ref{eq1})
\begin{equation}
\mathbf{P}\left(  \sup_{s\leq t}X_{s}^{(a)}>b\right)  \leq\exp\left(  \frac
{q}{1+ 1/p } \right)  \left(  \frac{t B(a) \left( p!\right) ^{1/p}}{ab}\right)
^{q/(1+1/p)}. \label{inq2}%
\end{equation}
Let $\beta< \alpha$ be defined later. If both $a(t)x$ and $t$ are small
enough, we get from (\ref{eq:B-asy}), (\ref{eq:a-prop1}), and (\ref{beta})
that for some $K > 0$
\[
\frac{t B(a) \left(  p!\right)  ^{1/p}}{ab} \leq K x^{\beta- \alpha}.
\]
Substituting back into (\ref{inq2}) we obtain
\[
\mathbf{P}\left(  \sup_{s\leq t}X_{s}^{(a)}>b\right)  \leq\exp\left(  \frac
{q}{1+ 1/p } \right)  \left(  Kx^{-\alpha+\beta}\right) ^{q/\left(
1+1/p\right) },
\]
which, by choosing $\beta>0$ small enough, implies (\ref{prime}). \qed

\begin{corollary}
\label{corr:2} Assume (\ref{LL}). For any $\beta\in(0, \alpha)$ there exists
$t_{0} > 0$ and $D > 0$ such that if $\max\{a(t) x, a(t), t \}< t_{0}$ and $x
\geq1$, then for any $\tau> 0$
\begin{equation}
\label{eq:cor2-ineq}\mathbf{P} \left(  \sup_{0 \leq s \leq t} \left(
-X_{t}^{(a(t)x)} \right)  > \tau a(t) x \right)  \leq\exp\left(  -D\tau
^{2}x^{\alpha-\beta}\right) .
\end{equation}
In particular, for any $0\leq s\leq t$%
\begin{equation}
\mathbf{P}\left(  -X_{s}^{(a)}>\tau a\right)  \leq\exp\left(  -D\tau
^{2}x^{\alpha-\beta}\right)  . \label{albet}%
\end{equation}
\end{corollary}

{\noindent\textit{Proof.}} Set $a=a(t) x$, and $b=\tau a$, with $x\geq1$, and
$\tau>0$. If both $a(t)x$ and $t$ are small enough, we get from
(\ref{eq:B-asy}), (\ref{eq:a-prop1}), and (\ref{beta}) that for some $D > 0$
\[
\frac{b^{2}}{ 2 t B(a)} \geq D \tau^{2} x^{\alpha-\beta}.
\]
Thus, an application of inequality (\ref{eq2}) implies (\ref{eq:cor2-ineq}). \qed

Recall from (\ref{eq:X-repr}) that $X^{-}$ is the spectrally negative part of
the L\'evy process $X$.

\begin{corollary}
\label{corr:3} Assume (\ref{LL}) and (\ref{eq:lambda+-}). For $\alpha= 1$
further assume that $\int_{[-1,0)} -y \Lambda(\mathrm{d} y) < \infty$. For
every $0<\varepsilon<1$ there exist $t_{0}>0$ and $x_{0}\geq1$, such that for
all $0<t<t_{0}$ and $x>x_{0}$
\begin{equation}
\mathbf{P}\left(  \sup_{s\leq t}X_{s}^{-}>\frac{\varepsilon}{4}a(t)x\right)
\leq e^{-x\varepsilon/8}. \label{claim}%
\end{equation}
\end{corollary}

{\noindent\textit{Proof.}} 
First note that if $\int_{[-1,0)} -y \Lambda(\dd y) < \infty$, in particular if $\alpha \leq 1$, then $-X_t^{-}$ is a subordinator, therefore the probability in question is 0.

Assume that $\alpha> 1$. Since $X_{s}^{-}$ is
a spectrally negative L\'{e}vy process for any $0<a\leq1$
\begin{equation}
\mathbf{P}\left(  \sup_{s\leq t}X_{s}^{-}>\frac{\varepsilon}{4}a(t)x\right)
\leq\mathbf{P}\left(  \sup_{s\leq t}X_{s}^{(-a)}- t \mu_{-}\left(  a\right)
>\frac{\varepsilon}{4}a(t)x\right)  , \label{pp}%
\end{equation}
where%
\[
\mu_{-}\left(  a\right)  =\int_{[-1,-a)} y\Lambda(\mathrm{d}y).
\]
Now%
\[
-t \mu_{-}\left(  a\right)  =- t \int_{[-1,-a)} y \Lambda(\mathrm{d}y) = t
\left(  a\overline{\Lambda}_{-}(a)-\overline{\Lambda}_{-}(1)+\int_{a}%
^{1}\overline{\Lambda}_{-}(y)\mathrm{d}y \right)  ,
\]
which by (\ref{eq:lambda+-}) and (\ref{LL}) for all small enough $a>0$ is for
some $C_{\alpha}>0$
\[
\leq C_{\alpha}at\overline{\Lambda}_{+}(a).
\]
Similarly, we can verify that for some $D_{\alpha}>0$%
\[
t B\left(  -a\right)  =t\int_{[-a,{0})} y^{2}\Lambda(\mathrm{d}y)\leq
D_{\alpha}a^{2}t\overline{\Lambda}_{+}(a).
\]
Setting $a=a (t)$ we see by (\ref{eq:a-prop1}) that for all $t>0$ small enough
for some $c_{\alpha}>0$ and $d_{\alpha}>0$ both%
\[
-\int_{[-1,-a)}ty\Lambda(\mathrm{d}y)\leq c_{\alpha} a(t) \ \text{ and }
\ tB\left(  -a\right)  \leq d_{\alpha} a^{2} (t) .
\]
Choose $x$ so large so that $\frac{\varepsilon}{8}x>\max\left\{  c_{\alpha
},2d_{\alpha}\right\}  $. Thus by (\ref{pp})
\[
\mathbf{P}\left(  \sup_{s\leq t}X_{s}^{-}>\frac{\varepsilon}{4}a(t)x\right)
\leq\mathbf{P}\left(  \sup_{s\leq t}X_{s}^{(-a\left(  t\right)  )}%
>\frac{\varepsilon}{8}a(t)x\right)  ,
\]
which by inequality (\ref{eq2}) in the $\sup_{s\leq t}X_{s}^{(-a)}$ case with
$a=a(t)$ and $b=\varepsilon a(t)x/8$ is
\[
\leq\exp\left(  -\frac{b^{2}}{2tB\left(  -a\right)  }\right)  \leq\exp\left(
-\frac{\left(  \frac{\varepsilon}{8}x\right)  ^{2}}{2d_{\alpha}}\right)
<\exp\left(  -\frac{\varepsilon}{8}x\right)  .
\]
This gives (\ref{claim}), with $x_{0} = 8 \varepsilon^{-1}\max\{ c_{\alpha}, 2
d_{\alpha}\}$ and for $t_{0} > 0$ sufficiently small.
\hspace*{1pt} \hfill \qed

\subsection{Four auxiliary lemmas}

The i.i.d.~counterpart of the next result is due to Bertoin, Lemma 1,
\cite{Bertoin98}.

\begin{lemma}
\label{lemma:max-proc1} Let $\alpha\in(0,2)$, $\alpha\neq1$. Assume (\ref{LL})
and (\ref{eq:lambda+-}). For any $0<\varepsilon<1$ there exist $A>0$,
$t_{0}>0$, $x_{0}\geq1$, $\alpha^{\prime}>\alpha$, such that if $t<t_{0}$,
$x>x_{0}$, $a(t)x<1$, and for $\alpha<1$ additionally assume $a(t)x<t_{0}$,
then
\begin{equation}
\mathbf{P}\left(  m_{t}\leq a(t)x(1-\varepsilon),\overline{X}_{t}%
>a(t)x\right)  \leq Ax^{-\alpha^{\prime}}. \label{ineqq}%
\end{equation}

\end{lemma}

{\noindent\textit{Proof.}} \textbf{Step 1.} Assume that $X_{t}$ is spectrally
positive. Note that in this case $\gamma_{-} = 0$. For $a\in(0,1)$, $c>0$ we
have
\[%
\begin{split}
&  \left\{  m_{t}\leq a,\overline{X}_{t}>c\right\}  =\left\{  N((0,t]\times
(a,\infty))=0,\ \overline{X}_{t}>c\right\} \\
&  =\{N((0,t]\times(a,\infty))=0\} \cap\left\{  \sup_{s\leq t}\left(
\gamma_{+}s+\int_{0}^{s}\int_{y\leq a} y \widetilde{N} (\mathrm{d}%
u,\mathrm{d}y)- s\int_{(a,1]} y\Lambda(\mathrm{d}y) \right)  >c\right\}  ,
\end{split}
\]
where the latter two events are independent. Therefore
\begin{equation}%
\begin{split}
&  \mathbf{P}\left(  m_{t}\leq a,\overline{X}_{t}>c\right) \\
&  =\mathbf{P}(m_{t}\leq a)\,\mathbf{P}\left(  \sup_{s\leq t}\left(
\gamma_{+}s+ \int_{0}^{s} \int_{(0,a]} y \widetilde{N}(\mathrm{d}%
u,\mathrm{d}y)- s \int_{(a,1]} y \Lambda(\mathrm{d}y)\right)  >c\right) \\
&  \leq\mathbf{P}\left(  \sup_{s\leq t}\left(  \gamma_{+}s+\int_{0}^{s}%
\int_{(0,a]} y \widetilde{N}(\mathrm{d}u,\mathrm{d}y)- s \int_{{(a,1]} } y
\Lambda(\mathrm{d}y)\right)  >c\right)  .
\end{split}
\label{eq:maxproc-ineq1}%
\end{equation}
Recall the definition from (\ref{NA}). We see by (\ref{eq:maxproc-ineq1}) that
when $1<\alpha<2$
\begin{equation*}
\mathbf{P}\left(  m_{t}\leq a,\overline{X}_{t}>c\right)  \leq\mathbf{P}\left(
\sup_{0\leq s\leq t}X_{s}^{(a)}>c\right) , \label{case1}%
\end{equation*}
as $\gamma_{+} = 0$ by (\ref{eq:gamma}), and when $0<\alpha<1$, again by
(\ref{eq:gamma})
\begin{equation*}
\mathbf{P}\left(  m_{t}\leq a,\overline{X}_{t}>c\right)  \leq\mathbf{P}\left(
\sup_{0\leq s\leq t}X_{s}^{(a)}>c-t\int_{(0,a]} y \Lambda\left(  \mathrm{d}
y\right)  \right)  . \label{case2}%
\end{equation*}
Fix $0<\varepsilon<1$. Put
\begin{equation}
a=a(t)x(1-\varepsilon) \quad\text{and } \ c=a(t)x. \label{eq:ab1}%
\end{equation}
In the case $0<\alpha<1$, by Karamata's theorem for $a>0$ small enough
\[
t \int_{(0,a]} y \Lambda(\mathrm{d}y) \leq t \int_{0}^{a} \overline{\Lambda
}_{+}(y) \mathrm{d}y \sim\frac{1}{1-\alpha} a t \overline{\Lambda}_{+}(a).
\]
Thus, by using (\ref{eq:a-prop1}) and the Potter bounds there exists a
$c_{1}>0$, such that for all $a=a(t)x(1-\varepsilon)>0$ small enough
\[
t \int_{(0,a]} y \Lambda(\mathrm{d}y) \leq a(t)x(1-\varepsilon) t
\overline{\Lambda}_{+} (a(t)x(1-\varepsilon)) \leq c_{1} a(t) x^{1-\alpha
+\alpha/2}.
\]
Hence for all $t>0$ small enough and $x$ large enough%
\begin{equation*}
a(t)x-t\int_{0}^{a}y\Lambda(\mathrm{d}y) >a(t)x\left(  1-c_{1}x^{-\alpha
+\alpha/2}\right)  > a(t) x \left(  1-\frac{\varepsilon}{2} \right) .
\label{a22}%
\end{equation*}
Therefore, for any $0<\alpha<2$, $\alpha\neq1$, \ there exist $t_{0}>0$,
$x_{0}\geq1$, $\alpha^{\prime}>\alpha$, such that for $0<t<t_{0}$, $x>x_{0}$,
and for $0<\alpha<1$ additionally assume $0<a(t)x<t_{0}$, we have with $a$ as
in (\ref{eq:ab1})
\begin{equation*}
\label{a23}\mathbf{P}\left(  m_{t}\leq a(t)x(1-\varepsilon),\overline{X}%
_{t}>a(t)x\right)  \leq\mathbf{P}\left(  \sup_{0\leq s\leq t}X_{s}^{(a)}>a(t)x
\left(  1-\frac{\varepsilon}{2} \right)  \right)  ,
\end{equation*}
which by inequality (\ref{prime}) is less than or equal to $Ax^{-\alpha
^{\prime}}$ for some $\alpha^{\prime}>\alpha$ and constant $A>0.$ This proves
(\ref{ineqq}).\medskip

\noindent\textbf{Step 2.} Finally, we extend the statement from spectrally
positive processes. Recall from (\ref{eq:X-repr}) that $X_{t}^{-}$ is the
spectrally negative part of $X_{t}$. Notice that by arguing as in Step 1, for
$a\in(0,1)$, $c>0 $ we have%
\[
\mathbf{P}\left(  m_{t}\leq a,\overline{X}_{t}>c\right)  \leq\mathbf{P}\left(
\sup_{s\leq t}\left(  X_{s}^{-}+\gamma_{+}s+X_{s}^{(a)}- s\int_{(a,1]} y
\Lambda(\mathrm{d}y)\right)  >c\right)
\]
In the case $0<\alpha<1$%
\[
0<-\int_{[-1,0)}y\Lambda(\mathrm{d}y)<\infty,
\]
so $-X_{t}^{-}$ is a subordinator and thus $X_{t}^{-}<0$ for any $t>0$.
Therefore, the result follows immediately from the $0<\alpha<1$ case of Step
1, since
\[%
\begin{split}
&  \sup_{s\leq t}\left(  X_{s}^{-}+\gamma_{+}s+X_{s}^{(a)}- s \int_{(a,1]}
y\Lambda(\mathrm{d}y)\right)  \leq\sup_{s \leq t}X_{s}^{(a)} + t \int_{(0,a]}
y\Lambda(\mathrm{d}y).
\end{split}
\]
On the other hand $\gamma_{+}=0$ in the case $1<\alpha<2$, thus
\[
\begin{split}
&  \left\{  \sup_{s\leq t}\left(  X_{s}^{-}+\gamma_{+}s+X_{s}^{(a)}%
-s\int_{(a,1]} y \Lambda(\mathrm{d}y) \right)  >a(t)x \right\} \\
&  \subset\left\{  \sup_{s\leq t}\left(  X_{s}^{(a)} - s\int_{(a,1]} y
\Lambda(\mathrm{d}y)\right)  > a(t) x \left(  1 - \frac{\varepsilon}{2}
\right)  \right\}  \cup\left\{  \sup_{s\leq t}X_{s}^{-}>\frac{\varepsilon}%
{2}a(t)x\right\}  .
\end{split}
\]
By a slight modification of first part of the proof given in Step 1, the
probability of the first event on the right-hand side of the last inclusion is
bounded by $Ax^{-\alpha^{\prime}}$ for some $\alpha^{\prime}>\alpha$ and
$A>0$, while the probability of the second event is exponentially small by
(\ref{claim}). \hspace*{1pt}\hfill\qed \medskip

\begin{lemma}
\label{lemma:max-proc2} Let $\alpha\in(0,2)$, $\alpha\neq1$. Assume (\ref{LL})
and (\ref{eq:lambda+-}). For any $0<\varepsilon<1$ there exist a constant
$A>0$, $t_{0}>0$, $x_{0}\geq1$, $\alpha^{\prime}>\alpha$, such that if
$t<t_{0}$, $x>x_{0}$, and $a(t)x<t_{0}$, then
\begin{equation}
\mathbf{P}\left(  \overline{X}_{t}\leq a(t)x(1-\varepsilon), m_{t}
>a(t)x\right)  \leq Ax^{-\alpha^{\prime}}. \label{alphaprime}%
\end{equation}

\end{lemma}

\noindent\textit{Proof.} \textbf{Step 1.} Assume first that $X$ is spectrally
positive. For $a\in(0,1)$ let $\tau=\tau_{a}=\inf\{s:\Delta X_{s}>a\}$. Then
$\{m_{t}>a\}=\{\tau\leq t\}$, and $\tau$ is exponentially distributed with
parameter $\overline{\Lambda}(a)$. Conditioning on $\tau$, and using
Proposition 0.5.2 in \cite{Bertoin}, for $b>0$%
\begin{equation}%
\begin{split}
&  \mathbf{P}\left(  \overline{X}_{t}\leq b,m_{t}>a\right)  = \mathbf{P}
\left(  \overline{X}_{t}\leq b,\tau\leq t\right) \\
&  \leq\mathbf{P}\left(  X_{\tau}\leq b,\tau\leq t\right) \\
&  =\int_{0}^{t} \mathbf{P}\left(  \gamma_{+} s + \int_{0}^{s} \int_{(0,a]}
y\widetilde{N} (\mathrm{d}u,\mathrm{d}y) - s \int_{(a,1]} y\Lambda
(\mathrm{d}y) \leq b-a \right)  \overline{\Lambda}_{+}(a) e^{-\overline
{\Lambda}_{+}(a)s}\mathrm{d}s\\
&  \leq\left(  1-e^{-\overline{\Lambda}_{+}(a)t}\right)  \sup_{s\leq
t}\mathbf{P}\left(  X_{s}^{(a)}+\gamma_{+}s- s \int_{(a,1]} y\Lambda
(\mathrm{d}y)\leq b-a\right)  .
\end{split}
\label{eq:cond}%
\end{equation}
Put
\begin{equation}
b=a(t)x(1-\varepsilon) \quad\text{and }\ a=a(t)x. \label{eq:ab3}%
\end{equation}
For $\alpha\in(1,2)$ integration by parts and Karamata's theorem give
\[%
\begin{split}
\int_{(z,1]} y\Lambda(\mathrm{d}y)  &  =\int_{z}^{1}\overline{\Lambda}%
_{+}(y)\mathrm{d}y+z\overline{\Lambda}_{+}(z)-\overline{\Lambda}_{+}(1)\\
&  \sim\frac{\alpha}{\alpha-1}z\overline{\Lambda}_{+}(z),\quad\text{ as
}z\downarrow0.
\end{split}
\]
Moreover, by (\ref{eq:a-prop1}) and by Potter's bounds, there exist $t_{0}>0$
and $x_{0}\geq1$ such that for $t<t_{0}$, $x>x_{0}$, $a(t)x<t_{0}$
\[
t\overline{\Lambda}_{+}(a)a=t\overline{\Lambda}_{+}(a(t))\frac{\overline
{\Lambda}_{+}(a)}{\overline{\Lambda}_{+}(a(t))}a\leq2x^{-\alpha/2}a(t)x.
\]
Therefore for any $0<\varepsilon<1$ fixed, there exist $t_{0}>0$ and
$x_{0}\geq1$ such that for $0\leq s\leq t<t_{0}$, $x>x_{0}$, $a(t)x<t_{0}$,
\begin{equation*}
b-a+s\int_{(a,1]} y \Lambda(\mathrm{d}y) \leq-\frac{\varepsilon}{2}a(t)x.
\label{eq:a-b1}%
\end{equation*}
For $\alpha\in(0,1)$ by the definition of $\gamma_{+}$ in (\ref{eq:gamma}),
simply $X_{s}^{(a)}+\gamma_{+}s - s\int_{(a,1]} y\Lambda(\mathrm{d} y) \geq
X_{s}^{(a)}$.

Summarizing, for any $\alpha\in(0,1) \cup(1,2)$, $0\leq s\leq t<t_{0}$ and
$x>x_{0}$, we get the bound
\[
\mathbf{P}\left(  X_{s}^{(a)}+\gamma_{+}s- s\int_{(a,1]} y\Lambda(\mathrm{d}
y)\leq b-a\right)  \leq\mathbf{P}\left(  X_{s}^{(a)}\leq-\frac{\varepsilon}%
{2}a(t)x\right)  .
\]
Inequality (\ref{albet}) gives for any choice of $0<\beta<\alpha$ there exist
$t_{0} > 0$, such that for $0<t< t_{0}$ and $0<a (t) x <t_{0}$ with $x\geq1$
\begin{equation*}
\mathbf{P}\left(  -X_{s}^{(a\left(  t\right)  x)}>\frac{\varepsilon}%
{2}a(t)x\right)  \leq\exp\left(  -D\left(  \frac{\varepsilon}{2}\right)
^{2}x^{\alpha-\beta}\right)  =:\exp\left(  -D\left(  \frac{\varepsilon}%
{2}\right)  ^{2}x^{\delta}\right)  , \label{AB}%
\end{equation*}
which is clearly stronger than (\ref{alphaprime}).\smallskip

\noindent\textbf{Step 2.} We extend the proof to the general case. As in
(\ref{eq:cond})%
\begin{equation}%
\begin{split}
&  \mathbf{P}\left(  \overline{X}_{t}\leq b,m_{t}>a\right) \\
&  \leq\left(  1-e^{-t\overline{\Lambda}_{+}(a)}\right)  \sup_{s\leq
t}\mathbf{P}\left(  X_{s}^{-}+X_{s}^{(a)}+\gamma_{+}s-s\int_{a}^{1}%
y\Lambda(\mathrm{d}y)\leq b-a\right) \\
&  \leq\left(  1-e^{-t\overline{\Lambda}_{+}(a)}\right)  \sup_{s\leq t}
\bigg[\mathbf{P}\left(  X_{s}^{(a)}+\gamma_{+}s- s\int_{(a,1]} y\Lambda
(\mathrm{d}y) \leq\frac{b-a}{2} \right)  + \mathbf{P}\left(  X_{s}^{-}%
\leq\frac{b-a}{2} \right)  \bigg].
\end{split}
\label{eq:cond2}%
\end{equation}
The first term in the square bracket is exponentially small by the first part
of the proof, where $a,b$ are as in (\ref{eq:ab3}).

Note that $-X_{s}^{-}$ in the second term is a spectrally positive L\'{e}vy
process, therefore we can use the methods of the first part of the proof of
Lemma \ref{lemma:max-proc1}. Let $m_{t}^{-}=\sup_{s\leq t}|X_{s}^{-}%
-X_{s-}^{-}|$, we have
\begin{equation}%
\begin{split}
&  \mathbf{P}\left(  -X_{s}^{-}> \frac{a-b}{2} \right)  =\mathbf{P}\left(
-X_{s}^{-}> \frac{\varepsilon}{2} a\left(  t\right)  x\right) \\
&  \leq\mathbf{P}\left(  m_{s}^{-}>\frac{\varepsilon}{4} a\left(  t\right)
x\right)  +\mathbf{P}\left(  -X_{s}^{-}>\frac{\varepsilon}{2} a\left(
t\right)  x, \ m_{s}^{-} \leq\frac{\varepsilon}{4} a\left(  t\right)  x
\right)  .\label{two}%
\end{split}
\end{equation}
For the first term in (\ref{two}) we have by (\ref{eq:lambda+-}),
(\ref{eq:a-prop1}), and (\ref{kappa})
\[
\begin{split}
\mathbf{P}\left(  m_{s}^{-}>\frac{\varepsilon}{4}a(t)x\right)   &  =
1-\exp\left\{  -s\Lambda_{-}(a(t)x\varepsilon/4)\right\} \\
&  \leq c_{1} t \Lambda_{+} (a(t)x\varepsilon/8) \leq c_{2} x^{-\alpha7/8},
\end{split}
\]
whenever $a(t) x$ and $t$ are small enough. For the second term in
(\ref{two}), by assumption (\ref{eq:lambda+-}), Lemma \ref{lemma:max-proc1} is
applicable, therefore it is of order $x^{-\alpha^{\prime}}$ for some
$\alpha^{\prime}> \alpha$. Finally, note that the first factor in the
right-hand side of (\ref{eq:cond2})
\[
1-e^{-t\overline{\Lambda}_{+}(a)} \leq2 t\overline{\Lambda}_{+}(a) \leq c_{3}
x^{-7 \alpha/8},
\]
and the result follows. \hspace*{1pt} \hfill\qed\medskip

The next result is the continuous analogue of Lemma 2, \cite{Bertoin98}.
Recall the notation in (\ref{eq:M-Y-def}).

\begin{lemma}
\label{lemma:MX} Let $\alpha\in(0,2)$, $\alpha\neq1$. Assume (\ref{LL}) and
(\ref{eq:lambda+-}). For any $y> 0$ and $\delta\in(0,1)$
\[
\lim_{u \downarrow0} \mathbf{P} \left(  M_{u} \leq(- \log u )^{1/\alpha} y (
1- \delta), Y_{u} > (- \log u )^{1/\alpha} y \right)  = 0.
\]

\end{lemma}

%\begin{proof}
{\noindent\textit{Proof.}} For $0<q<1$ consider the sequence $q^{n}$. We have
for $s \in[q^{n+1},q^{n}]$
\[
\frac{m_{s}}{a(s)}\geq\frac{m_{q^{n+1} }}{a(q^{n+1} )}\frac{a(q^{n+1}
)}{a(q^{n} )}.
\]
As $a$ is regularly varying, the second factor in the lower bound converges to
$q^{1/\alpha}$. Therefore for any $q^{1/\alpha}>\varepsilon_{1}>0$ there is a
$t_{0}>0$ such that for all $0< u \leq u^{\prime}\leq t_{0}$
\begin{equation}
\sup_{u \leq s\leq u^{\prime}}\frac{m_{s}}{a(s)} \geq\max_{n^{\prime}\leq
n\leq n_{u}} (q^{1/\alpha}-\varepsilon_{1}) \frac{m_{q^{n} }}{a(q^{n}
)}=:(q^{1/\alpha}-\varepsilon_{1}) \widetilde{M}_{u}, \label{m1}%
\end{equation}
where
\begin{equation}
\label{n}n_{u}=\left\lceil \log u/\log q\right\rceil \quad\text{and} \quad
n^{\prime}= n_{u^{\prime}} = \left\lceil \log u^{\prime}/\log q\right\rceil .
\end{equation}
Since $\overline{X}_{u}$ is monotone increasing, for $s\in[q^{n+1} ,q^{n} ]$
\[
\frac{\overline{X}_{s}}{a(s)}\leq\frac{\overline{X}_{q^{n} }}{a(q^{n} )}%
\frac{a(q^{n} )}{a(q^{n+1} )}.
\]
Similarly, the second factor in the upper bound converges to $q^{-1/\alpha}$.
Thus for any $q^{1/\alpha}>\varepsilon_{1}>0$ there is a $t_{0}>0$ such that
for all $0< u \leq u^{\prime}\leq t_{0}$%
\begin{equation}
\sup_{u \leq s\leq u^{\prime}} \frac{\overline{X}_{s}} {a(s)} \leq
(q^{-1/\alpha}+\varepsilon_{1}) \max_{n^{\prime} -1 \leq n\leq n_{u}}
\frac{\overline{X}_{q^{n} }}{a(q^{n})} =: (q^{-1/\alpha}+\varepsilon
_{1})\widetilde{Y}_{u}, \label{v1}%
\end{equation}
with $n^{\prime}$ and $n_{u}$ as above. Note that for any $0<u^{\prime}<1$
fixed
\begin{equation}
\sup_{u^{\prime}\leq s\leq1}\frac{m_{s}}{a(s)}=O_{\mathbf{P}}(1),\quad
\sup_{u^{\prime}\leq s\leq1}\frac{\overline{X}_{s}}{a(s)}=O_{\mathbf{P}}(1).
\label{eq:mV-sb}%
\end{equation}
Keeping in mind that $0<\delta<1$ is fixed, we can choose $q<1$ close to $1$
and $\varepsilon_{1}<q^{1/\alpha}$ small so that
\begin{equation}
(1-\delta)(q^{1/\alpha}-\varepsilon_{1})^{-1}<(q^{-1/\alpha}+\varepsilon
_{1})^{-1}. \label{eq:qeps}%
\end{equation}
Then choose $t_{0}$ such that both (\ref{m1}) and (\ref{v1}) hold true. This
choice will permit us to use Lemma \ref{lemma:max-proc1}. We see for $0 < u
\leq u^{\prime}\leq t_{0}$
\[%
\begin{split}
&  \mathbf{P}\left(  M_{u}\leq(- \log u )^{1/\alpha}y(1-\delta),Y_{u}>(- \log
u )^{1/\alpha}y\right) \\
&  \leq\mathbf{P}\left(  M_{u}\leq(- \log u )^{1/\alpha}y(1-\delta
),\sup_{u\leq s\leq u^{\prime}} \frac{\overline{X}_{s}}{a(s)}>(- \log
u)^{1/\alpha}y\right) \\
&  \qquad+\mathbf{P}\left(  \sup_{u^{\prime}\leq s\leq1}\frac{\overline{X}%
_{s}}{a(s)}>(- \log u )^{1/\alpha}y\right)
\end{split}
\]
which by (\ref{eq:mV-sb}), as $u\downarrow0$,
\[%
\begin{split}
&  =\mathbf{P}\left(  M_{u} \leq(- \log u )^{1/\alpha}y(1-\delta), \sup_{u
\leq s\leq u^{\prime}} \frac{\overline{X}_{s}}{a(s)}> (- \log u )^{1/\alpha
}y\right)  +o(1)\\
&  \leq\mathbf{P}\left(  \sup_{u \leq s\leq u^{\prime}}\frac{m_{s}}{a(s)}%
\leq(- \log u )^{1/\alpha}y(1-\delta),\sup_{u \leq s\leq u^{\prime}}%
\frac{\overline{X}_{s}}{a(s)}>(- \log u )^{1/\alpha}y\right)  +o(1)\\
&  \leq\mathbf{P}\left(  \widetilde{M}_{u}\leq\frac{(- \log u )^{1/\alpha
}y(1-\delta)}{q^{1/\alpha}-\varepsilon_{1}},\widetilde{Y}_{u}>\frac{(- \log u
)^{1/\alpha}y}{q^{-1/\alpha}+\varepsilon_{1}}\right)  +o(1),
\end{split}
\]
where at the last inequality we used (\ref{m1}) and (\ref{v1}).

We apply Lemma \ref{lemma:max-proc1} with
\begin{equation}
\label{eq:def_txeps}t=q^{n} ,\ x=x(u)=\frac{(- \log u )^{1/\alpha}%
y}{q^{-1/\alpha}+\varepsilon_{1}},\ \varepsilon= 1-\frac{q^{-1/\alpha
}+\varepsilon_{1}}{q^{1/\alpha}-\varepsilon_{1}}(1-\delta),
\end{equation}
and note that $\varepsilon\in(0,1)$ by (\ref{eq:qeps}). By Lemma
\ref{lemma:max-proc1} there exist $0<t_{1}\leq t_{0}$, $x_{1}>0$, and
$\alpha^{\prime}>\alpha$ such that if $t<t_{1}$, $x>x_{1}$, $a(t)x<t_{1}$
then
\[
\mathbf{P}\left(  \frac{m_{q^{n} }}{a(q^{n} )}\leq\frac{(- \log u )^{1/\alpha
}y(1-\delta)}{q^{1/\alpha}-\varepsilon_{1}},\frac{\overline{X}_{q^{n} }%
}{a(q^{n} )}>\frac{(- \log u )^{1/\alpha}y}{q^{-1/\alpha}+\varepsilon_{1}%
}\right)  \leq2 x^{-\alpha^{\prime}}.
\]
With $x(u)$ in (\ref{eq:def_txeps}), for $u>0$ define
\begin{equation}
\eta=\eta(u)=\min\{n:\,a(q^{n} )x(u)<t_{1}\}. \label{eta2}%
\end{equation}
Using Potter's bounds for $z$ small enough $a(z)\leq z^{1/(2\alpha)}$, thus
\[
a(q^{n} )x=a(q^{n} )\,\frac{(- \log u )^{1/\alpha}y}{q^{-1/\alpha}+\delta_{1}%
}\leq\frac{y}{y^{-1/\alpha}+\varepsilon_{1}}q^{n/(2\alpha)}(-\log
u)^{1/\alpha}.
\]
For $n \geq4(\log q^{-1})^{-1}\,\log\log u^{-1}$ we have
\[
q^{n/(2\alpha)}(- \log u)^{1/\alpha}\leq(-\log u)^{-1/\alpha}%
\]
which tends to 0 as $u \downarrow0$. Therefore, for $u$ small enough
$\eta(u)\leq4(\log q^{-1})^{-1}\log\log u^{-1}$. Recall $n_{u}$ and
$n^{\prime}$ from (\ref{n}). Simply,
\[%
\begin{split}
&  \mathbf{P}\left(  \widetilde{M}_{u} \leq\frac{(- \log u )^{1/\alpha
}y(1-\delta)}{q^{1/\alpha}-\varepsilon_{1}}, \widetilde{Y}_{u}>\frac{(- \log u
)^{1/\alpha}y}{q^{-1/\alpha}+ \varepsilon_{1}}\right) \\
&  \leq\mathbf{P}\left(  \overline{X}_{q^{n^{\prime}}}>t_{1}\right) \\
&  \qquad+\sum_{n=\eta(u)}^{n_{u}}\mathbf{P}\left(  \frac{m_{q^{n} }}{a_{q^{n}
}}\leq\frac{(- \log u )^{1/\alpha}y(1-\delta)}{q^{1/\alpha}-\varepsilon_{1}%
},\frac{\overline{X}_{q^{n} }}{a_{q^{n} }}>\frac{(- \log u )^{1/\alpha}%
y}{q^{-1/\alpha}+\varepsilon_{1}}\right) \\
&  \leq\mathbf{P}\left(  \overline{X}_{q^{n^{\prime}}}>t_{1}\right)  +2n_{u}
x^{-\alpha^{\prime}},
\end{split}
\]
where the second term goes to 0 as $u \downarrow0$ for any $u^{\prime}$. To
finish the proof note that for $t_{1}>0$ fixed as $u^{\prime}\downarrow0$
(thus $n^{\prime}\rightarrow\infty$) we have $\mathbf{P}(\overline
{X}_{q^{n^{\prime}}}>t_{1})\rightarrow0$.
%\end{proof}
\hspace*{1pt} \hfill\qed

\begin{lemma}
\label{lemma:MX2} Let $\alpha\in(0,2)$, $\alpha\neq1$. Assume (\ref{LL}) and
(\ref{eq:lambda+-}). For any $y>0$ and $0<\delta<1$
\[
\lim_{u \downarrow0}\mathbf{P}(Y_{t}\leq(- \log u )^{1/\alpha}y (1-\delta),
M_{u}>(- \log u )^{1/\alpha}y)=0.
\]
\end{lemma}

{\noindent\textit{Proof.}} The proof follows the steps of the previous proof,
so we only sketch it.

For $0<q<1$ consider the sequence $q^{n}$. For any $q^{1/\alpha}%
>\varepsilon_{1}>0$ there is a $t_{0}>0$ such that for all $0< u \leq
u^{\prime}\leq t_{0}$
\begin{equation}
\sup_{u\leq s\leq u^{\prime}}\frac{m_{s}}{a(s)}\leq\max_{n^{\prime} -1 \leq
n\leq n_{u}}(q^{-1/\alpha}+\varepsilon_{1})\frac{m_{q^{n}}}{a(q^{n}%
)}=:(q^{-1/\alpha}+\varepsilon_{1})\widetilde{M}_{u}, \label{m2}%
\end{equation}
and
\begin{equation}
\sup_{t\leq s\leq t^{\prime}}\frac{\overline{X}_{s}}{a(s)}\geq\max_{n^{\prime
}\leq n\leq n_{t}}(q^{1/\alpha}-\varepsilon_{1})\frac{\overline{X}_{q^{n}}%
}{a(q^{n})}=:(q^{1/\alpha}-\varepsilon_{1})\widetilde{Y}_{u}, \label{v2}%
\end{equation}
where $n^{\prime}$ and $n_{u}$ are defined as in (\ref{n}).

Choose $q<1$ close to $1$ and $\varepsilon_{1}<q^{1/\alpha}$ so small that
(\ref{eq:qeps}) holds. Then choose $t_{0}$ such that both (\ref{m2}) and
(\ref{v2}) hold true. This choice will permit us to use Lemma
\ref{lemma:max-proc2}. We see for $0<u \leq u^{\prime} \leq t_{0}$
\[%
\begin{split}
&  \mathbf{P}\left(  Y_{u} \leq(- \log u )^{1/\alpha}y(1-\delta),M_{u}>(- \log
u )^{1/\alpha}y\right) \\
&  \leq\mathbf{P}\left(  Y_{u}\leq(- \log u )^{1/\alpha}y(1-\delta
),\sup_{u\leq s\leq u^{\prime}}\frac{m_{s}}{a(s)}>(- \log u )^{1/\alpha
}y\right) \\
&  \qquad+\mathbf{P}\left(  \sup_{u^{\prime}\leq s\leq1}\frac{m_{s}}{a(s)}>(-
\log u )^{1/\alpha}y\right) ,
\end{split}
\]
where the second term goes to 0 by (\ref{eq:mV-sb}). For the first term by
(\ref{m2}) and (\ref{v2}) we have
\[%
\begin{split}
&  \mathbf{P}\left(  Y_{u}\leq(- \log u )^{1/\alpha}y(1-\delta),\sup_{u\leq
s\leq u^{\prime}}\frac{m_{s}}{a(s)}>(- \log u )^{1/\alpha}y\right) \\
&  \leq\mathbf{P}\left(  \widetilde{Y}_{u}\leq\frac{(- \log u )^{1/\alpha
}y(1-\delta)}{q^{1/\alpha}-\varepsilon_{1}},\widetilde{M}_{u}>\frac{(- \log u
)^{1/\alpha}y}{q^{-1/\alpha}+\varepsilon_{1}}\right)  .
\end{split}
\]
Choose $t, x, \varepsilon$ as in (\ref{eq:def_txeps}). Using Lemma
\ref{lemma:max-proc2} we can show there exist $A>0$, $\alpha^{\prime}>\alpha$,
$0<t_{1}\leq t_{0}$, $x_{1}\geq1$, and $A>0$ such that if $t<t_{1}$, $x>x_{1}%
$, $a(t)x<t_{1}$ then
\[
\mathbf{P}\left(  \frac{\overline{X}_{q^{n}}}{a(q^{n})}\leq\frac{(- \log u
)^{1/\alpha}y(1-\varepsilon)}{q^{1/\alpha}-\varepsilon_{1}},\frac{m_{q^{n}}%
}{a(q^{n})}>\frac{(- \log u )^{1/\alpha}y}{q^{-1/\alpha}+\varepsilon_{1}%
}\right)  \leq Ax^{-\alpha^{\prime}}.
\]
For $\eta(u)$ as in (\ref{eta2}), as in the previous proof for $u$ small
enough $\eta(u) \leq4 ( \log q^{-1})^{-1}\log\log u^{-1}$.We obtain
\[%
\begin{split}
&  \mathbf{P}\left(  \widetilde{Y}_{u}\leq\frac{(- \log u )^{1/\alpha
}y(1-\delta)}{q^{1/\alpha}-\varepsilon_{1}},\widetilde{M}_{u}>\frac{(- \log u
)^{1/\alpha}y}{q^{-1/\alpha}+\varepsilon_{1}}\right) \\
&  \leq\mathbf{P}\left(  m_{q^{n^{\prime}}}>t_{1}\right)  +\sum_{n=\eta
(u)}^{n_{u}}\mathbf{P}\left(  \frac{\overline{X}_{q^{n}}}{a(q^{n})}\leq
\frac{(- \log u )^{1/\alpha}y(1-\varepsilon)}{q^{1/\alpha}-\varepsilon_{1}%
},\frac{m_{q^{n}}}{a(q^{n})}>\frac{(- \log u )^{1/\alpha}y}{q^{-1/\alpha
}+\varepsilon_{1}}\right) \\
&  \leq\mathbf{P}\left(  m_{q^{n^{\prime}}}>t_{1}\right)  +A n_{u}%
x^{-\alpha^{\prime}},
\end{split}
\]
where the second term goes to $0$ as $u\downarrow0$ for any $u^{\prime}$. To
finish the proof note that for $t_{1}>0$ fixed as $t^{\prime}\downarrow0 $
(thus $n^{\prime}\rightarrow\infty$) we have $\mathbf{P}(m_{q^{n^{\prime}}%
}>t_{1})\rightarrow0$. \hspace*{1pt} \hfill\qed \medskip

\subsection{Proof of Theorem \ref{thm:Y-lim}}

Now we are ready to prove Theorem \ref{thm:Y-lim}. Let $0<\varepsilon<1$ be
arbitrary. Simply,
\[%
\begin{split}
&  \mathbf{P}(M_{t}\leq(- \log t )^{1/\alpha}x(1-\varepsilon))\\
&  =\mathbf{P}(M_{t}\leq(- \log t )^{1/\alpha}x(1-\varepsilon),Y_{t}\leq(-
\log t )^{1/\alpha}x)\\
&  \quad+\mathbf{P}(M_{t}\leq(- \log t )^{1/\alpha}x(1-\varepsilon),Y_{t}>(-
\log t )^{1/\alpha}x).
\end{split}
\]
By Theorem \ref{thm:M-lim} the left-hand side converges to $\exp
\{-[(1-\varepsilon)x]^{-\alpha}\}$, and by Lemma \ref{lemma:MX} the second
term in the right-hand side tends to 0. Therefore
\[
\lim_{t\downarrow0}\mathbf{P}(M_{t}\leq(- \log t )^{1/\alpha}x(1-\varepsilon
),Y_{t}\leq(- \log t )^{1/\alpha}x)=\exp\left\{  -[(1-\varepsilon)x]^{-\alpha
}\right\}  ,
\]
thus
\begin{equation}
\liminf_{t\downarrow0}\mathbf{P}(Y_{t}\leq(- \log t )^{1/\alpha}x)\geq
e^{-x^{-\alpha}}. \label{eq:Yliminf}%
\end{equation}

On the other hand, for $0<\varepsilon<1$
\[%
\begin{split}
&  \mathbf{P}\left(  Y_{t}\leq(- \log t )^{1/\alpha}x(1-\varepsilon)\right) \\
&  =\mathbf{P}\left(  Y_{t}\leq(- \log t )^{1/\alpha}x(1-\varepsilon),M_{t}>(-
\log t )^{1/\alpha}x\right) \\
&  \qquad+\mathbf{P}\left(  Y_{t}\leq(- \log t )^{1/\alpha}x(1-\varepsilon
),M_{t}\leq(- \log t )^{1/\alpha}x\right)  .
\end{split}
\]
Here the first term on the right-hand side goes to 0 by Lemma \ref{lemma:MX2},
and by Theorem \ref{thm:M-lim}
\[
\limsup_{t\downarrow0}\mathbf{P}\left(  Y_{t}\leq(- \log t )^{1/\alpha
}x(1-\varepsilon),M_{t}\leq(- \log t )^{1/\alpha}x\right)  \leq e^{-x^{-\alpha
}}.
\]
Combining this with (\ref{eq:Yliminf}) the result follows.
%\end{proof}
%\hspace*{1pt} \hfill\qed

\subsection{Proof of Theorem \ref{thm:Y-lima1}}

In the $\alpha= 1$ case the result follows similarly, only a minor change is
needed in the proof, because one cannot choose the centering to be zero. Note
that Theorem \ref{thm:M-lim}, Proposition \ref{prop}, and Corollaries
\ref{corr:1}, \ref{corr:2}, and \ref{corr:3} hold for any $\alpha\in(0,2)$. 
Recalling the definition of the centering in (\ref{eq:c-a1}),
introduce the notation 
\begin{equation} \label{eq:hatX}
\widehat X_t = \sup_{s \leq t} (X_s - c(s)), \quad t \geq 0.
\end{equation}
Lemma \ref{lemma:max-proc1} remains true in the following form.

\begin{lemma}
\label{lemma:max-proc1-a1} Assume (\ref{LL}) with $\alpha= 1$,
(\ref{eq:lambda+-}), and $\int_{[-1,0)} -y \Lambda(\mathrm{d} y) < \infty$.
For any $0<\varepsilon<1$ there exist $A>0$, $t_{0}>0$, $x_{0}\geq1$,
$\alpha^{\prime}>1 $, such that if $t<t_{0}$, $x>x_{0}$, $a(t)x<t_{0}$, then
\[
\mathbf{P}\left(  m_{t}\leq a(t)x(1-\varepsilon),
\widehat X_t > a(t)x \right)  \leq Ax^{-\alpha^{\prime}}.
\]
\end{lemma}

{\noindent\textit{Proof.}} 
\textbf{Step 1.}
First let $X_{t}$ be spectrally positive. Note that in this case $\gamma_- = 0$. 
For $a = a(t) x (1 -\varepsilon) \in (0,1)$, $c = a(t) x$ we have
\[%
\begin{split}
&  \left\{  m_{t}\leq a, \widehat X_t >c\right\}  
=\left\{  N((0,t]\times (a,\infty))=0,\ \widehat X_t >c\right\} \\
&  =\{N((0,t]\times(a,\infty))=0\}
\cap\left\{  \sup_{s\leq t} \left(  \gamma_{+} s+
\int_{0}^{s} \int_{y\leq a} y \widetilde{N} (\mathrm{d}u,\mathrm{d}y)-
s \int_{(a,1]} y\Lambda(\mathrm{d}y) 
- c(s) \right)  >c\right\}  ,
\end{split}
\]
where the latter two events are independent. Therefore
\begin{equation}%
\mathbf{P}\left(  m_{t}\leq a,\widehat{X}_{t}>c\right) 
\leq \mathbf{P} \left(  \sup_{s\leq t}\left(  \gamma_{+} s+\int_{0}^{s}%
\int_{(0,a]} y \widetilde{N}(\mathrm{d}u,\mathrm{d}y)-
s \int_{{(a,1]} } y \Lambda(\mathrm{d}y) - c(s) \right)  >c\right)  .
\label{eq:maxproca1-ineq1}%
\end{equation}

Recall the definition of $\gamma_+$ and the centering in (\ref{eq:gamma}) and in (\ref{eq:c-a1}).
Since $a(t) x(1-\varepsilon) > a(s)$, for $s \leq t$ and $x$ large enough, 
if $\int_{(0,1]} y \Lambda( \dd y) = \infty$
we obtain
\begin{equation} \label{eq:a1-ac=1}
\begin{split}
\gamma_+ s - s \int_{{(a,1]} } y \Lambda(\mathrm{d}y) - c(s) 
& = - s \int_{{(a,1]} } y \Lambda(\mathrm{d}y)  + 
s \int_{(a(s), 1]} y \Lambda(\dd y) \\
& = s \int_{(a(s), a(t) x (1-\varepsilon)]} y 
\Lambda(\dd y) > 0.
\end{split}
\end{equation}
While, if $\int_{(0,1]} y \Lambda(\dd y) < \infty$, 
\begin{equation} \label{eq:a1-ac=2}
\begin{split}
\gamma_+ s - s \int_{{(a,1]} } y \Lambda(\mathrm{d}y) - c(s) 
& =
s \int_{(0,1]} y \Lambda(\dd y) - s \int_{{(a,1]} } y \Lambda(\mathrm{d}y)  - 
s \int_{(0,a(s)]} y \Lambda(\dd y) \\
& = s \int_{(a(s), a(t) x (1-\varepsilon)]} y 
\Lambda(\dd y) > 0.
\end{split}
\end{equation}
Therefore in both cases we get the same term. Next, we claim that
\begin{equation} \label{eq:a1-ineqac1}
\sup_{s \leq t}  
s \int_{(a(s), a(t) x]} y \Lambda(\dd y) \leq 
\frac{\varepsilon}{2} a(t) x.
\end{equation}
We have for $x > 1$ and $t \geq s  > 0$ small
\begin{equation} \label{eq:aux2-a1}
\begin{split}
s \int_{(a(s), a(t)x]}  y \Lambda (\dd y) 
& = s \left( \int_{a(s)}^{a(t)x} \overline \Lambda_+(y) \dd y
- \ell(a(t)x) + \ell(a(s)) \right) \\
& \leq s \int_1^{\frac{x a(t)}{a(s)}}  
\frac{\ell(a(s) u)}{u} \dd u + s \ell(a(s) ) .
\end{split}
\end{equation}
By Potter's bounds, whenever $a(t) x$ is small enough
\[
\int_1^{\frac{x a(t)}{a(s)}} 
\frac{\ell(a(s) u)}{\ell(a(s)) u} \dd u \leq 
\int_1^{\frac{x a(t)}{a(s)}} 2 {u}^{-1/2} \dd u 
< 4 \sqrt{x} \, \sqrt{ \frac{a(t)}{a(s)}}.
\]
Substituting back into (\ref{eq:aux2-a1}) 
and using that $\overline \Lambda_+(a(t)) = \ell(a(t)) / a(t) \sim t^{-1}$ by (\ref{eq:a-prop1}), we obtain uniformly in $s \leq t$
\[
\begin{split}
& s \int_{(a(s), a(t)x]}  y \Lambda (\dd y) 
\leq s \ell(a(s)) \left( 
4  \sqrt{x} \sqrt{\frac{a(t)}{a(s)}} + 1 \right) \\
& \leq 5 x^{-1/2} x a(t) \leq \frac{\varepsilon}{2} x a(t)
\end{split}
\]
for $x$ large enough and $t$ small enough.
This proves (\ref{eq:a1-ineqac1}).

Using the bound (\ref{eq:a1-ineqac1}) in inequality (\ref{eq:maxproca1-ineq1}) we obtain 
\[
\mathbf{P}\left(  m_{t}\leq a,\widehat{X}_{t}>c\right) 
\leq \mathbf{P} \left(  
\sup_{s \leq t} X_s^{(a)} > a(t) x \left( 1 - \frac{\varepsilon}{2} \right) \right),
\]
and the result follows from (\ref{prime}).

\noindent \textbf{Step 2.} The extension to the general case is immediate now, because $- X^{-}_t$ is a subordinator by our assumption
$\int_{[-1,0)} - y \Lambda(\dd y) < \infty$.
\hspace*{1pt} \hfill\qed

The corresponding version of Lemma \ref{lemma:max-proc2}
also holds. Recall the definition in (\ref{eq:hatX}).

\begin{lemma}
\label{lemma:max-proc2-a1} Assume (\ref{LL}) with $\alpha= 1$ and
(\ref{eq:lambda+-}). For any $0<\varepsilon<1$ there exist a constant $A>0$,
$t_{0}>0$, $x_{0}\geq1$, $\alpha^{\prime}>1$, such that if $t<t_{0}$,
$x>x_{0}$, and $a(t)x<t_{0}$, then
\[
\mathbf{P}\left(  
\widehat X_t \leq a(t)x(1-\varepsilon),
 m_{t} >a(t)x\right)  \leq Ax^{-\alpha^{\prime}}.
\]
\end{lemma}

\noindent\textit{Proof.} 
%\textbf{Step 1.} 
Assume first that $X_t$ is spectrally
positive. For $a\in(0,1)$ let $\tau=\tau_{a}=\inf\{s:\Delta X_{s}>a\}$.
As in the proof of Lemma \ref{lemma:max-proc2} for $b>0$
\begin{equation*}%
\begin{split}
&  \mathbf{P}\left(  \widehat{X}_{t}\leq b,m_{t}>a\right)  = 
\mathbf{P} \left(  \widehat{X}_{t}\leq b,\tau\leq t\right) \\
&  \leq\mathbf{P} \left(  X_{\tau} - c(\tau) \leq b,\tau\leq t\right) \\
&  =\int_{0}^{t} 
\mathbf{P}\left(  \gamma_{+} s + \int_{0}^{s} \int_{(0,a]}
y\widetilde{N} (\mathrm{d}u,\mathrm{d}y) - s \int_{(a,1]} y\Lambda
(\mathrm{d}y) - c(s) \leq b-a \right)  \overline{\Lambda}_{+}(a) e^{-\overline
{\Lambda}_{+}(a)s}\mathrm{d}s\\
&  \leq\left(  1-e^{-\overline{\Lambda}_{+}(a)t}\right)  \sup_{s\leq
t}\mathbf{P}\left(  X_{s}^{(a)}+\gamma_{+}s- s \int_{(a,1]} y\Lambda
(\mathrm{d}y) - c(s) \leq b-a\right)  .
\end{split}
\label{eq:conda1}%
\end{equation*}
Put $b=a(t)x(1-\varepsilon)$ and $a=a(t)x$.
From (\ref{eq:a1-ac=1}) and (\ref{eq:a1-ac=2}) we obtain
\[
\mathbf{P} \left(  X_{s}^{(a)}+
\gamma_{+} s - s\int_{(a,1]} y\Lambda(\mathrm{d} y) - c(s) \leq b-a\right)  \leq\mathbf{P}\left(  X_{s}^{(a)}\leq b-a\right)  .
\]
Therefore, the result follows as in the proof of Lemma \ref{lemma:max-proc2}.\smallskip

The general case follows exactly as in the proof of Lemma \ref{lemma:max-proc2}.
%\noindent\textbf{Step 2.}  
\hspace*{1pt} \hfill\qed \medskip

After having the appropriate versions of Lemma \ref{lemma:max-proc1} and
\ref{lemma:max-proc2} the proof of the theorem is identical to the proof in
the $\alpha\neq1$ case.

\subsection{Proof of Theorem \ref{thm:DE}}

We shall prove that
\begin{equation}
\lim_{t\downarrow0}\mathbf{P}\left(  Y_{t}=\sup_{t\leq s\leq1}\frac
{X_{s}-c(s)}{a(s)}\right)  =1,\label{onez}%
\end{equation}
which clearly implies the theorem. 

First assume that $\alpha\neq1$, in which case $c(s)=0$. Note that
\[
Y_{t}=\sup_{t\leq s\leq1}\frac{\sup_{u\leq s}X_{u}}{a(s)}\geq\sup_{t\leq
s\leq1}\frac{X_{s}}{a(s)}=:Z_{t}.
\]
Assume that $Y_{t}>Z_{t}$. Then $Y_{t}=\frac{X_{u_{0}}}{a(s_{0})}$, for some
$s_{0}\in\lbrack t,1]$ and $u_{0}\leq s_{0}$. Since $Y_{t}>Z_{t}$, we have
$u_{0}<t$, thus the monotonicity of $a$ implies $Y_{t}=\overline{X}_{t}/a(t)$.
Therefore
\[
\mathbf{P}(Y_{t}>Z_{t})\leq\mathbf{P}\left(  Y_{t}=\frac{\overline{X}_{t}%
}{a(t)}\right).
\]
Now for all $t>0$ and $x>0$
\begin{equation} \label{eq:pt}
\mathbf{P}\left(  Y_{t}=\frac{\overline{X}_{t}}{a(t)}\right)  \leq
\mathbf{P}\left(  Y_{t}\leq x(-\log t)^{1/\alpha}\right)  +\mathbf{P}\left(
Y_{t}=\frac{\overline{X}_{t}}{a(t)},\frac{\overline{X}_{t}}{a(t)}\geq x(-\log
t)^{1/\alpha}\right)  =:p_{t}\left(  x\right)  .
\end{equation}
By Theorem 2 for all $x>0$ the first term on the right-hand side tends to
$\exp\left(  -x^{-\alpha}\right)$, which converges to $0$ as $x \downarrow 0$.
Next we show that $\overline{X}_{t}/a(t)$ is stochastically bounded.
By (\ref{eq:X-repr})
\[
\frac{1}{a(t)} \sup_{s \leq t} X_s \leq 
\frac{1}{a(t)} \sup_{s \leq t} X^+_s 
+ \frac{1}{a(t)} \sup_{s \leq t} X^{-}_s.
\]
The second term is stochastically bounded by (\ref{claim}), while the first term is stochastically bounded since the process%
\[
\frac{X^+_{ts}}{a(t)}, \ 0\leq s\leq1,
\]
converges weakly in $D_{0}[ 0, 1]$. (See Remark \textit{(iv)}
on page 322 of Maller and Mason \cite{MM08} and the methods of the proofs of
Proposition 4.1 and Corollary 4.2 of Maller and Mason \cite{MM10}.) Thus the
second term in (\ref{eq:pt}) converges to $0$ for all $x>0$. We see now that
\[
\limsup_{t\downarrow0}\mathbf{P}\left(  Y_{t}=\frac{\overline{X}_{t}}%
{a(t)}\right)  \leq\lim_{x\downarrow0}\limsup_{t\downarrow0}p_{t}\left(
x\right)  =\lim_{x\downarrow0}\exp\left(  -x^{-\alpha}\right)  =0,
\]
which implies that
\[
\lim_{t\downarrow0}\mathbf{P}\left(  Y_{t}=Z_{t}\right)  =1,
\]
which is (\ref{onez}). 
\smallskip

For $\alpha = 1$ the proof is almost identical. There is a small difference in the stochastic boundedness of $\widehat X_t/a(t)$. Note that
\[
\frac{1}{a(t)} \sup_{s \leq t} (X_s - c(s) ) \leq 
\frac{1}{a(t)}  \sup_{s \leq t} (X_s^+ - c(s) ) 
+ \frac{1}{a(t)}  \sup_{s \leq t} X_s^{-}.
\]
The second term is again stochastically bounded by (\ref{claim}), while for the first it follows from the convergence $(X_t - c(t)) /a(t)$ as above.
\hspace*{1pt} \hfill$\square$

\bigskip

\noindent\textbf{Acknowledgment.} P\'{e}ter Kevei is supported by the
J\'{a}nos Bolyai Research Scholarship of the Hungarian Academy of Sciences, by
the NKFIH grant FK124141, and by the EU-funded Hungarian grant
EFOP-3.6.1-16-2016-00008. David Mason's visit at the Bolyai Institute was
supported by the Ministry of Human Capacities, Hungary grant
20391-3/2018/FEKUSTRAT. David Mason thanks the hospitality of the Bolyai
Institute. He also acknowledges the support of a 2019 UDARF Research Fund Grant.

%\bibliographystyle{abbrv}
%\bibliography{DE}
%\end{document}

\end{document}